\documentclass[11pt]{article} 
\usepackage{amscd} 
\usepackage{amssymb} 
\usepackage{amsmath} 
\usepackage{amsthm} 
\usepackage{amsbsy} 
\usepackage{graphicx} 
\usepackage[all]{xy}

\input{psfig.tex} 
\input{psfig.sty} 
 
\newcommand{\R}{\ensuremath{\mathbf R}} 
 
\newcommand{\Z}{\ensuremath{\mathbf Z}} 
 
\newcommand{\Pro}{\ensuremath{\mathbf{RP}}} 
\newcommand{\zd}{\ensuremath{\Z/2\Z}} 
 
\newcommand{\zo}{\ensuremath{\Z/8\Z}} 
\newcommand{\su}{\ensuremath{S^1}}

\newcommand{\tchi}{\ensuremath{\tilde{\chi}}} 
\newcommand{\Mk}{\ensuremath{M_k}} 
 
\newtheorem{teo}{Theorem}[section] 
\newtheorem{prop}[teo]{Proposition} 
\newtheorem{coroll}[teo]{Corollary} 
\newtheorem{definiz}[teo]{Definition} 
\newtheorem{oss}[teo]{Remark} 
\newtheorem{lem}[teo]{Lemma} 
\newtheorem{ex}[teo]{Example}

\newenvironment{dimostraz}{\noindent{\em Proof.}\ }{\vspace{0.5cm}}

\def \2{{\bf Z}/2{\bf Z}} 
\def \o{{\mathcal O}} 
\newcommand{\QP}{\ensuremath{Q{\mathbf{RP}}^\infty}} 
 
\begin{document} 
\title{The graded cobordism group of codimension-one immersions} 
\author{ 
\begin{tabular}{cc} 
 Louis Funar\thanks{Partially supported by a Canon grant 1999-2000.} 
 & Rosa Gini\\ 
\small \em Institut Fourier BP 74, UMR 5582 CNRS, &\small \em Department of 
Mathematics, \\ \small \em University of Grenoble I, &\small \em 
University of Pisa,\\ \small \em 38402 Saint-Martin-d'H\`eres cedex, France 
&\small \em via Buonarroti 2, 56127 Pisa, Italy\\ \small \em e-mail: {\tt 
funar@fourier.ujf-grenoble.fr} & \small \em e-mail: {\tt 
gini@dm.unipi.it} \\ 
\end{tabular} 
} 
\date{} 
 
\maketitle 
{\abstract 
The cobordism group $N(M^n)$ of codimension-one immersions in the $n$-manifold
$M^n$ has a natural filtration induced by any cellular decomposition.
The problem addressed in this paper
is the explicit computation of the graded group $gr^*N(M^n)$.
We introduce some new invariants for immersions  
enlightening the Atiyah-Hirzebruch spectral sequence associated to $N(M)$,
which are of combinatorial-geometric nature. Explicit computations are
developed for $n\leq7$, and the group structure is also investigated for
orientable 4-manifolds.

}

\vspace{0.3cm} 
 
\vspace{0.3cm} 
 
\tableofcontents 
 
\section{Introduction}

The classification of manifold immersions in codimension at least one 
up to cobordism was reduced to a homotopy problem by the results of 
\cite{VogCI,WelCGI}. These techniques are however awkward to apply if one 
wants to get effective results. The classification up to regular 
homotopy is also a homotopy problem which is closely related to the 
previous one. For instance two immersed surfaces in ${\bf R^3}$ of 
the same topological type are regularly homotopic if and only if they 
are cobordant (see \cite{PinRHC}). This subject received recently more 
attention (see e.g. \cite{NowQPR}). 
 
The group  $P_n$ of codimension-one immersions in the $n$-sphere 
up to cobordism  is the 
 the $n$-th stable  homotopy group of ${\bf 
  RP}^{\infty}$ and they were computed by Liulevicius (\cite{LiuTHA}) 
for $n\leq 9$. Explicit classifications for regular homotopy 
equivalence of immersed  surfaces in 3-manifolds were first obtained 
by Hass and Hughes in \cite{HaHIST} and Pinkall (see \cite{PinRHC}). 
 
The cobordism group $N(M^3)$ of immersed surfaces in the 3-manifold $M^3$ 
was computed  geometrically  by Benedetti and Silhol 
(\cite{BeSSPS}). Let $M^3$ be a compact oriented 3-manifold 
and  $f$ a smooth codimension-one immersion of a (compact) 
surface $F^2$ in $M^3$. Fixing a Spin 
structure on $M^3$ one has a Pin structure induced on $F^2$ which 
defines a ${\bf Z}/4{\bf Z}$-valued quadratic form on $H_1(F^2, 
{\bf Z}/2{\bf Z})$, by counting how the immersion $f$ twists the 
regular neighborhoods of 1-cycles in $F^2$. 
There is then an isomorphism between $N(M^3)$ and 
$H_1(M^3, {\bf Z}/2{\bf Z}) \oplus H_2(M^3, {\bf Z}/2{\bf Z}) \oplus 
{\bf Z}/8{\bf Z}$, the last being endowed with a twisted product. 
The isomorphism sends an immersion $f$ into the triple consisting of 
the homology class of the double points locus, the homology class of 
the image of $f$, and the Arf invariant of  the quadratic form from 
above. A similar result holds for nonorientable 3-manifolds $M^3$, but 
the factor ${\bf Z}/8{\bf Z}$ is replaced now by ${\bf Z}/2{\bf Z}$ 
(see \cite{rosa}). Notice that the factor ${\bf Z}/8{\bf Z}$ is 
nothing but $N(S^3)$ so the two results above can be stated in a unitary 
way by considering $H^3(M^3, {\bf Z}/8{\bf Z})$. 
 
One would like to have a similar description for the group $N(M^n)$ 
in all dimensions $n$. More motivation for that is the result of 
\cite{FunCIM}, which relates the cobordism group $N(M^n)$ to the set 
$CB(M)$ of cubulations of the manifold $M^n$ modulo a set of 
combinatorial moves analogous to Pachner's move on 
simplicial complexes. 
We refer to \cite{FunCIM} for an extensive 
discussion of this problem, due to Habegger 
(see problem 5.13 from \cite{KirPLT}). 
 
One remarks first that there is a natural grading $gr^*$ on $N(M^n)$ 
induced by 
a cellular decomposition. The Atiyah-Hirzebruch spectral sequence (see 
\cite{HilGCT}) has its second term $E_2^{p,q}=H^q(M^n, P_p)$ and converges 
to the graded $N(M^n)$. However one has only very few informations about the 
differentials in this sequence, hence the direct use of this approach fails. 
One develops then a combinatorial way to settle this question. 
 
One obtains thus the extension of  the computations 
of the graded group $gr^*N(M^n)$ up to dimension 
$4$, and up to dimension 7 under a mild homological condition. 
However the techniques one uses are 
different from those of Benedetti and Silhol, 
though as they are still geometric in nature.

The main theorem is the following: 
 
\begin{teo}\label{main} Let $M$ be a closed $n$-manifold, $n\leq7$. Then 
\[ 
gr^*(N(M))=H^1(M,P_1)\times EH^2(M)\times H^3(M,P_3)\times\dots\times 
H^n(M,P_n). 
\] 
\end{teo} 
 
\noindent The subgroup $EH^2(M)\subseteq H^2(M,P_2)$ is defined in 
section~\ref{EH} 
and is computed for any 
$n\leq4$, and in some cases for $n\in\{5,6,7\}$, as follows: 
\[ 
EH^2(M)=\left\{ 
\begin{array}{ll} 
H^2(M,P_2)&\text{ if }n\leq3 \text{ or $n=4$ and $M$ non-orientable}\\ 
\{x\in H^2(M,P_2); \, x\cup x =0\}&\text{ if $n=4$ and $M$ orientable or $n>4$ and condition $(*)$ holds} 
\end{array} 
\right. 
\] 
where the condition $(*)$ for $n$-manifolds  with $n\in\{5,6,7\}$ is 
that $M$ is orientable and 
${\rm Ext}(H_3(M),\zo)=0$. 
 
The whole theory is explicitly developed for closed manifolds. 
However the present methods can be applied to non-compact manifolds simply by 
substituting the ordinary cohomology with the cohomology with compact 
support.

\vspace{0.2cm} 
 
\noindent{\bf A sketch of proof}. We briefly summarize the guiding line 
of the paper. %
We introduce in section~\ref{filtrationsec} the natural filtration of $N(M)$ that 
gives rise to the graded group, and prove that $F^k$ can be interpreted 
geometrically as the subgroup of $N(M)$ of those immersions avoiding 
the $k$-skeleton up to cobordism. This property is independent of the 
cellular decomposition. In section~\ref{cohomological_invariants} we 
define in a geometric way an injective homomorphism 
\[ 
\widetilde{\chi}^k:F^k/F^{k+1}\longrightarrow H^k(M,P_k)/NEH^k(M), 
\] 
where $NEH^k(M)\subseteq H^k(M,P_k)$ is a subgroup of cocycles related 
to some particular null-cobordant immersions. We then introduce an 
obstruction theory that permits to study the inverse map to 
$\widetilde{\chi}^k$, in particular to determine its image. In 
section~\ref{EH} the theory is applied to explicit computations that 
provide the image of $\widetilde{\chi}^k$ 
up to $n=7$, and in section~\ref{NEH} it is proven that $NEH^k(M)=0$ 
for all $k$ in all $n$-manifold $M$ up to $n=7$. These computations 
prove the main theorem.

\vspace{0.2cm} 
 
\noindent{\bf Acknowledgements}.  We are indebted to Riccardo Benedetti, 
Octav Cornea, Takuji Kashiwabara, Valentin Po\'enaru, 
John Scott-Carter and Pierre Vogel 
for valuable discussions and suggestions.  Part of 
this work has been done when the first author visited Tokyo Institute 
of Technology, which he wishes to thank for the support and hospitality, 
and especially to Teruaki Kitano and Tomoyoshi Yoshida. 
 
\section{The groups $P_n$ and $Q_n$ in low dimension}\label{PnQn} 
 
Let $M$ be a $n$-dimensional  manifold. Consider the set of 
immersions $f:F\rightarrow M$ with $F$ a closed $(n-1)$-manifold. 
Impose on it the following relation: $(F,f)$ is {\em cobordant} to 
$(F',f')$ if there exist 
 a cobordism $X$ between $F$ and $F'$, that is, a 
compact $n$-manifold $X$ with boundary $F\sqcup F'$, and an immersion $\Phi:X\rightarrow M\times I$, transverse to 
the boundary, 
 such that $\Phi_{|F}=f\times\{0\}$ 
and $\Phi_{|F'}=f'\times\{1\}$. 
 
Once the manifold $M$ is fixed, the 
set $N(M)$ of cobordism classes of codimension-one immersions  in $M$ is 
an abelian 
group with the composition law given by disjoint union. 
 
In this paper we mainly deal with the cobordism group of immersions in 
manifolds of dimension less or equal than 7. This is due to the 
fact that the groups $P_n:=N(S^n)=N(\R^n)$, which are always finite, 2-torsion groups 
(see \cite{WelCGI}), are particularly simple for 
$n<7$, as is shown in the following table (see \cite{LiuTHA}): 
 
\vspace{0.5cm} 
\begin{center} 
\begin{tabular}{|c|c|c|c|c|c|c|c|c|c|}\hline 
$n$ &  1 & 2 & 3 & 4 & 5 & 6 & 7 & 8 & 9\\ \hline 
$P_n$ &$ {\bf Z}/2{\bf Z}$ & ${\bf Z}/2{\bf Z}$ & ${\bf Z}/8{\bf Z}$ 
&${\bf Z}/2{\bf Z} $ &$ 0$ & ${\bf Z}/2{\bf Z}$ & ${\bf Z}/16{\bf 
  Z}\oplus {\bf Z}/2{\bf Z} 
$ &$ ({\bf Z}/2{\bf Z})^{\oplus 3}$ & $({\bf Z}/2{\bf Z})^{\oplus 4} $\\ \hline 
\end{tabular} 
\end{center} 
\vspace{0.5cm} 
 
\noindent The simplest case is $n=5$, but the cases $n=1,2,4,6$ are also 
very  easy to handle. Indeed consider the classical invariant 
\[ 
\theta_n:P_n\longrightarrow \zd, 
\] 
that associates to an immersion the number of $n$-tuple points modulo 2, 
or equivalently, the homology class modulo 2 of the set of $n$-tuple 
points (as an element of $H_0(S^n,\zd)$). It is well-known that 
$\theta_n$ is an isomorphism for $n=1,2,4,6$ (see for example 
\cite{CarGBS}). The group $P_2$ is generated by the immersion 
8 which looks like the figure eight in the plane, while the group $P_4$ 
is generated by an immersion of $S^3$ with a single quadruple point. 
 
By definition of cobordism 
an immersion $f$ of a compact $(n-1)$-manifold in $S^n$ represents the trivial element of $P_n$ if 
and only if there exists an immersion of a compact $n$-manifold with boundary 
in $S^n\times I$ transverse to the boundary. It is clear that this is equivalent 
to ask that 
$f$ 
bounds in $D^{n+1}$. We generalize slightly this condition. 
 
\begin{definiz} Denote by 
$Q_n$ the subgroup of $P_n$ of those immersions in 
$S^n$ bounding an immersion of a compact $n$-manifold with boundary 
in a  $(n+1)$-manifold with boundary $S^n$. 
\end{definiz} 
 
The reason why $\theta_n$ being an isomorphism makes computations 
easier also when dealing with immersions in manifolds that are not spheres 
amounts to the following easy proposition, lying on the elementary 
but fundamental 
 fact that compact 1-manifolds with boundary have an even 
number of points as boundary. 
 
\begin{prop}\label{bounding} Let $f: F\longrightarrow S^n$ be a a codimension-one 
immersion, and let $n$ be such that either $\theta_n$ is an 
isomorphism or $P_n$ 
is trivial. 
Then $f$ bounds an immersion in a $(n+1)$-manifold whose boundary is $S^n$ 
if and only if it bounds in $D^{n+1}$, that is 
$Q_n=0$. 
\end{prop} 
 
\begin{dimostraz} 
If $n$ is such that $P_n=0$ there is nothing to prove. Assume then that $\theta_n$ is an 
isomorphism and 
let $f$ be a immersion in $S^n$. Suppose that there exists a $(n+1)$-manifold 
$N$ with boundary $S^n$ and a generic codimension-one immersion $g$ in $N$ transverse 
to the boundary and such that $g\cap{\partial N}=f$. 
The set of $n$-tuple 
points of $g$ is the immersion in $N$ of a compact 1-manifold with boundary, bounding the 
set of $n$-tuple points of $f$. But a 1-manifold has an even number of 
points as boundary, so $\theta_n(f)=0$, that is, $f$ represents 
the trivial element of $P_n$. 
~$\square$ 
 
\end{dimostraz} 
 
\noindent This proposition does not apply, for example, for $n=3$. The group 
$P_3=\zo$ is generated by the left Boy immersion, which has a single triple point. The even elements 
have no triple points. Canonical representatives for these classes are 
the immersions in $\R^3$ obtained by rotating an 8 on the $xz$ plane 
with the double point in $(1,0,0)$ around the $z$ axis, while rotating 
it in its own plane of half a twist, a whole twist and respectively 
three halves of 
twists (see \cite{PinRHC}). These immersions, whose Arf invariants are 
2,4 and 6 respectively, have a circle of double points. The immersion 
similarly constructed that makes no rotations is null-cobordant.

From 
proposition~\ref{bounding} it immediately follows: 
 
\begin{coroll} For any $n$ the group $Q_n$ is contained in 
$\ker\,\theta_n$. 
\end{coroll} 
 
\noindent In computing $Q_3$ we will make use first of a natural 
way of producing codimension-one immersions. 
 
\begin{definiz}\label{decorating} {\rm Let $N$ be a $n$-manifold 
    (possibly with boundary) and 
$S$ an embedded codimension-$k$ submanifold. Assume that the 
structure group of the 
normal bundle $\nu$ to $S$ in $M$  can be reduced to 
the group of symmetries of an element $f\in P_k$, 
and choose such a reduction denoted by $\mu$. 
There is then a canonical embedding of 
$f$ in each fiber $\R^k$ of $\nu$ giving rise to a  sub-fibration 
of $\nu$ with fiber $f$.  The total space of the last fibration 
is an immersion in the tubular neighborhood of $S$, hence in $N$, 
which will be called the {\em immersion obtained by decorating 
$S$ with $f$}\/ and will be denoted by $S\ltimes_\mu f$ or by $S\ltimes f$ if the reduction is not 
relevant to the context. 
The cobordism class of this immersion is clearly 
independent of the representative of the cobordism class of $f$, 
but might depend on $\mu$. 
} 
\end{definiz} 
 
\noindent For example remark that the symmetry group of the 8 in 
$\R^2$ is the same as the symmetry group of a line, so decorating 
a codimension-two submanifold with an 8 is the same as choosing a 
field of lines on the submanifold itself, if any. 
In particular the canonical representatives of the even 
elements of $P_3$ can be considered as $S^1\ltimes_{\mu_m} 8$, 
where $S^1$ is the standard circle in $\R^3$ and $\mu_m$ is the 
line field on the circle making $m$ halves of of twist, 
$m=0,1,2,3$, hence ${\rm Arf}(S^1\ltimes_{\mu_m} 8) =2m$. Remark 
that 
in an orientable $n$-manifold 
any simple curve has 
trivial normal bundle, hence it can be decorated by any element of 
$P_{n-1}$, the cobordism class of the resulting immersion only 
depending on the choice of the trivialization. On the other side 
in a non-orientable $n$-manifold a simple curve non-trivially intersecting the 
orientation cycle cannot be decorated by an element of $P_{n-1}$ 
not admitting a reflection in its symmetry group.

\begin{prop}\label{Q3} $Q_3=2P_3$. 
 
\end{prop} 
 
\begin{dimostraz} 
Let $M$ be the non-orientable $S^3$ bundle on \su, let $e^4$ be a 
4-ball in $M$ intersecting $\su$ and $N=M\setminus int(e^4)$. Let $\gamma$ be $\su\setminus 
(int(e^4)\cap\su)$. The normal bundle to $\gamma$ is trivial hence one can decorate it 
with left Boy immersions. Denote by $g$ the resulting immersion and by $f$ 
its intersection with the three sphere $\partial e^4$. Remark that when seen in the boundary of  $e^4$ (with any of its possible orientations) the two 
connected components of $f$ 
have the same orientation: otherwise the whole of \su\ could be 
decorated by left Boy immersions, which is not possible since the 
Boy immersions do not admit reflections in their symmetry group. 
Hence  $f$ determines a non-trivial element of 
$P_3$ that has Arf invariant 2 or 6, according to the orientation one 
has chosen on $\partial e^4$. But it is clear from its construction that $f$ 
bounds in $N=M\setminus int (e^4)$. So $f$ is in $Q_3$ and the same holds for the 
subgroup of $P_3$ it generates, that is, the whole of $2P_3$.~$\square$ 
 
\end{dimostraz} 
 
\noindent The proof of this proposition easily extends to the following 
statement: 
 
\begin{prop}\label{pari}  For any $n$ one has  $2P_n\subseteq Q_n$. 
Moreover any $f\in 2P_n$ 
bounds in any non-orientable $(n+1)$-manifold with boundary $S^n$.~$\square$ 
\end{prop} 
 
Set then  ${\cal P}_n=P_n/Q_n$ for the group of immersions in 
$S^{n}$ up to cobordism in  manifolds bounding 
two spheres. We proved the following: 
 
\vspace{0.5cm} 
\begin{center} 
\begin{tabular}{|c|c|c|c|c|c|c|}\hline 
$n$ &  1 & 2 & 3 & 4 & 5 & 6 \\ \hline 
${\cal P}_n$ &$ {\bf Z}/2{\bf Z}$ & ${\bf Z}/2{\bf Z}$ & ${\bf Z}/2{\bf Z}$ 
&${\bf Z}/2{\bf Z} $ &$ 0$ & ${\bf Z}/2{\bf Z}$\\ \hline 
\end{tabular} 
\end{center} 
\vspace{0.5cm}

\begin{oss}\label{bounding_in_orientable} {\rm  The immersion with invariant 4 
bounds in an orientable 4-manifold. Indeed,  let 
$L$ be a sphere corresponding to a complex line in ${\bf CP}^2$. 
There exists a normal  vector field $\mu$ on $L$ with a single zero. 
Let $e^4$ be a 4-ball  containing this zero. 
The restricted normal field $\tilde{\mu}$ then trivializes the normal bundle to $L\setminus (L\cap 
e^4)$ in $N={\bf CP}^2\setminus int(e^4)$, so one can define the immersion 
$g:=L\setminus (L\cap e^4)\ltimes_{\tilde{\mu}} 8$. A straightforward 
computation shows that $f:=g\cap\partial e^4$  is an immersion with Arf invariant 4, and 
clearly $f$ bounds in the orientable manifold $N$. 
} 
\end{oss} 
 
\section{A natural filtration of N(M)}\label{filtrationsec} 
 
The point of view from which we are able to tackle the computation of 
$N(M)$ is that of splitting it in pieces. At first sight the 
splitting depends on the cellular decomposition of $M$. 
 
\begin{definiz}{\rm Let $M$ be a $n$-manifold and let 
\[ 
M_0\subseteq\dots\subseteq M_k \subseteq\dots\subseteq M_n 
\] 
be a skeleton decomposition.  Let  $F^k\subset N(M)$, for $k\geq1$, be 
the set of immersions that up to 
cobordism do not intersect $M_{k-1}$. An immersion $f$ whose class 
belongs to $F^k$ will be said {\em $k$-admissible} if $f\cap 
M_{k-1}=\emptyset$. 
} 
 
\end{definiz} 
 
\noindent Remark that $F^k$ is a subgroup of $N(M)$ hence one has 
a filtration of $N(M)$. One has 
\[ 
N(M)=F^0\supseteq\dots\supseteq F^k \supseteq\dots\supseteq F^n, 
\] 
where $F^0$ was added for convenience of notation. 
This filtration comes in fact in a natural way from the 
algebraic-topological definition of $N(M)$. 
 
\begin{prop}\label{filtration} Under the Pontryagin-Thom construction we have 
\[ 
F^k=\{\varphi\in [M,\QP] \text{ such that }\varphi_{|M_{k-1}}\sim*\}=\ker 
([M,\QP]\rightarrow [M_{k-1},\QP]), 
\] 
where $\sim$ means homotopy and $*$ is the trivial based loop. 
\end{prop} 
 
\begin{dimostraz} If an immersion $f$ does not intersect $M_{k-1}$ then 
the map $\varphi_f$ associated by the Pontryagin-Thom construction is 
constant on $M_{k-1}$. On the other side, given a map $\varphi$ that is 
null-homotopic when restricted to $M_{k-1}$, 
consider 
the homotopy 
\[ 
H:M_{k-1}\times I\longrightarrow Q{\bf RP}^{\infty}, 
\] 
such that $H(-,0)=\varphi_{|M_{k-1}}$ and $H(-,1)$ is the constant map 
on $M_{k-1}$. 
The inclusion of the $(k-1)$-skeleton being a cofibration implies 
(see \cite{PicLHT}) that $H$ extends 
to $M\times I$.  Thus it gives a homotopy between $\varphi$ and a map 
$\varphi'$ defined on all of $M$,  whose restriction to $M_{k-1}$ is 
constant. Consider 
a closed regular neighborhood $N$ of $M_{k-1}$ 
in $M$. There exists then a global 
retraction $r:M\to M$ so that $r(N)=M_{k-1}$. 
 The map $\varphi''=\varphi'\circ r$ is thus constant on $N$ 
(in particular on $\partial N$) and 
is homotopic to $\varphi$. Therefore the Thom-Pontryagin 
construction  for the manifold with boundary 
$M-int(N)$ associates to the map $\varphi''_{|M\setminus int(N)}$ an 
immersion in $M\setminus int(N)$. When looking at that immersion as 
contained in $M$ it results as  a $k$-admissible 
representative of the class associated to $\varphi$, hence the claim is 
proved. 
 
We give a second proof introducing a recursive technique that will be 
often exploited in the sequel. We assume that the cellular 
decomposition is in particular a cubulation. 
Recall that a {\it cubical complex} is a complex  $K$ 
consisting of Euclidean cubes, 
such that the intersection of two cubes  is a finite union of cubes 
from  $K$, once a cube is in $K$ then all its faces belong to $K$, and no 
identifications of faces of the same cube are allowed. 
A {\it cubulation} of a manifold is specified by a cubical 
complex PL homeomorphic to the manifold.

Consider a class with a representative immersion $f$ for which 
the map $\varphi_f$ defined by the Pontryagin-Thom construction belongs 
to 
$\ker(N(M^n)\to [M_{k-1},\QP])$. We want to deform $f$ up to cobordism 
in such a way that the new representative does not intersect $M_{k-1}$. 
There exists a homotopy $H$ of $M$ so that 
the restriction of $\varphi_f\circ H|_{M_{k-1}}$ is 
null-homotopic. 
One uses now a recurrence on the degree $k$. 
If $k=1$ then it is obvious since the immersion can miss the 
$0$-skeleton by general position. Assume the claim is true for degree 
at most $k-1$. Then there exists a representative 
immersion so that $f\cap M_{k-2}^n=\emptyset$. 
This means that the intersection $f \cap e^{k-1}$ with any 
$(k-1)$-cell is a closed immersed submanifold lying in the interior 
of the cell. By hypothesis we can assume that 
${\varphi}_f|_{M_{k-1}}=1$, where 
$1$ denotes the constant (trivial) map. This means that 
$\varphi_{f\cap e^{k-1}}={\varphi}_f|_{e^{k-1}}=1$. 
The Thom-Pontryagin theory implies that the immersion 
$f\cap e^{k-1}$ is null cobordant. 
Consider a small regular neighborhood $V$ of $e^k$ in $M$, which is 
a product $e^k\times B^{n-k}$. Since $f \cap e^{k-1}$ is null 
cobordant there exists an immersion $g$ in $e^{k-1}\times I$ 
providing a null-cobordism for $f\cap e^{k-1}$. 
One uses $g \times \partial B^{n-k}$ to change the immersion 
in $V$ so that the new immersion misses $e^{k-1}$.~$\square$ 
\end{dimostraz}

Now it is  a classical result that such a filtration is independent on 
the cellular decomposition. Indeed  $[\cdot,\QP]$ is the 0-th degree of 
the  generalized cohomology theory associated to the suspension spectrum of 
$\QP$ 
\[ 
h^q(X)=\lim_{n\rightarrow\infty}[\Sigma^n 
X,\Sigma^{n+q}\Pro^{\infty}],\qquad  q\in\Z, 
\] 
and so the Atiyah-Hirzebruch spectral 
sequence converges to the graded group associated to the filtration of 
proposition~\ref{filtration}. 
The filtration being independent on the cellular 
decomposition for realizable generalized cohomology theories 
is  then  illustrated in the first and third chapter of the book of Hilton 
\cite{HilGCT}.

\noindent From now on we will often choose to use cubical 
decompositions of the 
manifold $M$. This will permit to perform recursive constructions in 
an easier way because a cubulation of $M$ induces in an obvious way 
one for  $M\times I$. 
 
\begin{oss} {\rm 
Most of the analysis that follows, and notably the definition of 
cohomological invariants and the application of obstruction 
theory, does not make use of any specific property of the generalized cohomology theory of cobordism 
groups of immersions, 
except perhaps the fact that it has finite coefficients. It is 
possible then that the same definitions apply usefully to other 
generalized cohomology theories. 
 
} 
\end{oss}

\section{Cohomological invariants}\label{cohomological_invariants}

Recall that a $k$-admissible immersion 
$f$ is such that $f\cap M_{k-1}=\emptyset$. In particular for any 
$k$-cell 
$e^k$ the intersection 
$f\cap e^k$ is contained in $int(e^k)$. 
If $e^k$ is oriented then  $f\cap e^k$ detects an element of $P_k$. 
One introduces  then the following 
geometric definition.  To any $k$-admissible immersion $f$ 
there is  associated a cochain $\chi_f^k\in C^k(M,P_k)$ the following 
way: 
\[ 
\chi^k_f(e^k):=f\cap e^k\in P_k. 
\] 
\begin{prop}For $k$-admissible $f$ the cochain $\chi^k_f$ is a cocycle. 
\end{prop} 
 
\begin{dimostraz} Given a $(k+1)$-cell $e^{k+1}$ recall 
that $\delta\chi^k_f(e^{k+1})=\chi^k_f(\partial 
e^{k+1})$ holds. Now denote 
$\partial e^{k+1}=\sum_{i\in {\cal I}(e^{k+1})}\varepsilon_ie_i^{k}$ 
where ${\cal I}(e^{k+1})$ is a finite set, $e_i^{k}$ are oriented 
$k$-cells (not necessarily different 
from each other) and $\varepsilon_i=\pm1$.  Consider the cell as a 
closed $(k+1)$-disk attached to $M_k$ by means of an attaching 
map, that results to be a homeomorphism when restricted to any 
connected component of the preimage of any $e_i^{k}$. 
Since $f\cap M_{k-1}$ is empty one can then pull 
back $f\cap {e}^{k+1}$ in the disk. 
The restriction of the resulting immersion to 
the boundary of the disk is then 
\[ 
\sum_{i\in {\cal I}(e^{k+1})}\varepsilon_i(f\cap e_i^{k})=\sum_{i\in {\cal I}(e^{k+1})} 
\varepsilon_i\chi^k_f(e_i^{k})\in N(S^k) 
\] 
and is trivial since it bounds in the disk. 
But this is $\chi^k_f(\partial 
e^{k+1})$, hence $\chi^k_f$ is a cocycle.~$\square$ 
\end{dimostraz} 
 
\noindent The argument of the previous proof will be repeatedly used. It is easy 
to visualize it when the cellular decomposition is a cubulation. For 
an immersion  $f$ and a $(k+1)$-cube $e^{k+1}$ the intersection $f\cap 
e^{k+1}$ is a cobordism to the empty set of $f\cap \partial e^{k+1}$. 
Thus the last one is trivial as an element of $P_{k}$.  When $f$ 
is $k$-admissible $f\cap \partial e^{k+1}$ splits as the sum (with 
signs) of $\chi^k_f(\text{faces})$. 
 
\begin{figure}[h] 
\begin{center} 
\includegraphics[width=4cm]{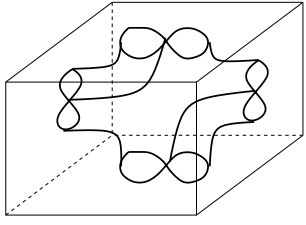} 
\caption{A $k$-admissible immersion restricted to a $(k+1)$-cube}\label{closed} 
\end{center} 
\end{figure}

Assume from now on that the cellular decomposition is  a cubulation. 
It is easy to see that if $f$ and $f'$ are $k$-admissible immersions 
admitting a cobordism $g$, which does not intersect $M_{k-1}\times I$, then 
$\chi^k_f=\chi^k_{f'}$. In fact for any $k$-cube $e^k$ the intersection 
$g\cap (e^k\times I)$ is a cobordism between $f\cap e^k$ and $f'\cap 
e^k$. 
The most natural question is whether $f$ 
and $f'$ admit a $k$-admissible cobordism, that is a cobordism 
not intersecting $M_{k-2}\times I$. This leads to the 
following result: 
 
\begin{prop}\label{condizione} Let $f$ and $f'$ be $k$-admissible cobordant immersions 
admitting a cobordism that is $k$-admissible in $M\times I$. Then 
$\chi^k_f$ and $\chi^k_{f'}$ are cohomologous. 
 
\end{prop} 
 
\begin{dimostraz} Assume $f$ and $f'$ transverse to the 
cubulation, and take a cobordism $g$ between them that is transverse 
to the standard cubulation of $M\times I$ associated 
with the cubulation of $M$. Our aim is to define a 
coboundary between $\chi^k_f$ and $\chi^k_{f'}$ by means of the 
same definition we used for $\chi^k$, but now applied to $g$. 
Remark that for any 
$k$-cube $e^k$ the intersection $g\cap(e^k\times I)$ is a cobordism to 
the empty set of $g\cap\partial (e^k\times I)$, thus the last one 
is then a trivial element  in $P_k$. 
Since  $g$ is 
$k$-admissible this element splits as the sum of $f\cap e^k-f'\cap e^k=\chi^k_f-\chi^k_{f'}$ 
and of $g\cap (\partial e\times I)$, see figure~\ref{closed}. 
We claim that the last 
summand is the coboundary of a $(k-1)$-cochain of $M$. 
 
Define a 
cochain $\psi$ in $C^{k-1}(M,P_k)$ this way. For any oriented 
$(k-1)$-cell $e^{k-1}$ consider the class of the immersion $g\cap 
(e^{k-1}\times I)$ where the orientation of $e^{k-1}\times I$ is such that it 
induces on $e^{k-1}$ the opposite of its orientation. With this 
convention  $g\cap(e^{k-1}\times I)$ is a well-defined 
element of $P_k$.

It is then easy to see that $\delta\psi (e^k)=g\cap(\partial e^k\times 
I)$, hence $\delta\psi=\chi^k_{f'}-\chi_{f}^k$.~$\square$ 
\end{dimostraz}

\noindent  It is 
immediate that two cobordant immersions which are 1-admissible 
have a 1-admissible cobordism between them.  For 
general $k$ we are not able to prove the analogous statement. 
We will face this problem 
gradually. If $g$ is a generic  cobordism between $k$-admissible 
immersions, and $e^s$ is an $s$-cube of $M$ we denote by $g(e^s)$ the 
immersion $g\cap(e^s\times I)$, where the orientation of $e^s\times I$ is such that it 
induces on $e^s$ the opposite of its orientation.  Now if $g$ is $s$-admissible then 
$g(e^s)$ is actually an immersion 
in the open $(s+1)$-disk $int(e^s)\times (0,1)$, hence 
$g(e^s)$ 
represents an element of $P_{s+1}$. Further if $g(e^{k-2})$ is empty for any 
$(k-2)$-cube then $g$ is $k$-admissible. 
 
\begin{prop}\label{oneadmissible} Let $f$ and $f'$ be 2-admissible cobordant immersions. 
Then there exists a 2-admissible cobordism between them. 
 
\end{prop} 
 
\begin{dimostraz} Let $g$ be a generic cobordism between $f$ and $f'$, 
transverse to the cubulation of $M\times I$ associated to a chosen 
cubulation of $M$. 
We want to prove that we can deform $g$ until $g(e^0)=\emptyset$ for 
any $e^0$. 
 
 Remark that $g(e^0)$  represents an element of $P_1$.  We 
claim that, up to modify $g$ if $M$ is compact, 
this element is trivial for any $e^0$. First remark that if 
$e^0$ and $f^0$ are two vertices of the cubulations that are connected 
by a edge $e^1$, then $g(e^0)$ and $g(f^0)$ are the same element of 
$P_1$. This is because $g(e^1)$ provides a cobordism between them and  this 
proves  in fact, $M$ being connected, that there is a well-defined 
$g^0\in P_1$ such that $g(e^0)=g^0$ for any vertex $e^0$ of $M$. Now 
suppose $M$ is not compact.  Since the domain of $g$ is compact 
its image cannot 
intersect all  edges of the type $e^0\times I$, hence $g^0$ is 
the trivial element of $P_1$, and the claim is proved in this case. If 
$M$ is compact, 
then $g^0$ might be non-trivial. But then consider a new cobordism 
obtained by adding a $M\times\{t\}$ to $g$.  Call $g$ the new 
cobordism and now $g^0$ is trivial as required. 
 
So we are ready to get rid of intersections of type $g(e^0)$. For any 
vertex $e^0$ there is a diffeomorphism of a neighborhood $U(e^0)\times 
I$ with $B^n\times I$ such that $g\cap D^n\times I$ is the inclusion of 
an even number of disks at levels $p_1,\dots p_{2n}$, since $g$ can be assumed to be transverse to $e^0\times I$. 
In this model cut 
the corresponding disks of radius $\frac{1}{2}$, and connect the holes in pairs 
by means of cylinders $\frac{1}{2}S^{n-1}\times [p_{2i-1},p_{2i}]$. The immersion 
obtained by repeating this construction in any vertex and then smoothing, that we still 
call $g$, satisfies $g(e^0)=\emptyset$ for any vertex $e^0$, hence $g$ 
is 2-admissible.~$\square$ 
\end{dimostraz} 
 
\noindent In general, given a $s$-admissible cobordism $g$ the elements $g(e^{s-1})$ 
for $s>1$ might be nontrivial. If they are trivial however it is 
possible to deform $g$ to a $(s+1)$-admissible cobordism. Consider the standard cubulation of $M\times I$ 
associated to a cubulation of $M$. We will say that cubes of the form 
$e\times I$ are {\em vertical}, cubes of the form $e\times \{0\}$ are {\em at the bottom} 
and cubes of the form $e\times \{1\}$ are {\em at the top}. In 
general cubes of the form $e\times \{t\}$ are {\em horizontal}. 
 
\begin{lem}\label{sadmissible} Let $f$ and and $f'$ be 
$(s+1)$-admissible immersions, let  $g$ be a cobordism between 
them 
that is $s$-admissible and such 
that for any $(s-1)$-cell $e^{s-1}$ of $M$ the immersion $g(e^{s-1})$ 
represents a trivial element of $P_{s}$. Then $g$ can be deformed to 
a $(s+1)$-admissible cobordism $g'$. 
 
If $g$ is a $s$-admissible cobordism between $f$ and $f'$, but 
$f$ alone is $(s+1)$-admissible, then $g$ can be  modified to 
a $(s+1)$-admissible cobordism $g'$ between $f$ and an immersion 
coinciding with $f'$ outside a neighborhood of $M_s$. 
 
%
%
%
 
\end{lem} 
 
\begin{dimostraz} By a construction analogous to that of proposition~\ref{oneadmissible} 
it is possible to get rid of intersections. 
We put again 
ourselves in a model, as follows: the normal bundle 
to $e^{s-1}$ in $M$ is trivial, 
take a trivialized neighborhood $U(e^{s-1})=e^{s-1}\times B^{n-s+1}$ 
and take its product with the interval $I$. From transversality one 
can suppose 
that $g\cap U(e^{s-1})\times I$ has the structure of a product 
$g(e^{s-1})\times B^{n-s+1}$. Consider a cobordism to the empty set of 
$g(e^{s-1})$, let $h$ 
be the embedding of this cobordism in $e^{s-1}\times \frac{1}{4} B^1 \times I\subset U(e^{s-1}) \times I$ 
and take the product 
$h\times \frac{1}{2} S^{n-s-1}$. Remark that the image of this product does not 
intersect $e^{s-1}\times I$, and that it intersects $g$ in $g(e^{s-1})\times 
\frac{1}{2} S^{n-s}$. Thus one can excise $g(e^{s-1})\times \frac{1}{2} B^{n-s+1}$ and 
glue back $h\times \frac{1}{2} S^{n-s}$. After repeating this construction 
for any $(s-1)$-cube and 
smoothing one gets a cobordism $g'$ between $f$ and $f'$ 
satisfying $g'(e^{s-1})=\emptyset$ for any 
$e^{s}$. 
 
If $f$ and $f'$ are both $(s+1)$-admissible then $g$ is clearly 
$(s+1)$-admissible, since the $s$-skeleton of $M\times I$ are the 
vertical $e^{s-1}\times I$ plus the two horizontal copies of 
$M_s$. 
 
If $f$ is $(s+1)$-admissible, but $f'$ intersects the $s$-skeleton, 
then we deform $g'$ further. Remark that for any $s$-cube $e^s$ of 
$M$ the cobordism $g'\cap (e^s\times I)$ is a cobordism to the 
empty set for $f'\cap e^s$. Hence this last immersion represents 
the trivial element of $P_s$. Then a surgery similar to the one of 
the first part of the proof leads to the excision of all of the 
intersections $g'\cap (e^s\times \{1\})$, and this proves the 
claim.~$\square$ 
 
\end{dimostraz} 
 
Since we cannot claim than any null-cobordant 
$k$-admissible immersion admits a $k$-admissible cobordism to the 
empty set one introduces the following definition. 
 
\begin{definiz}{\rm Set $NEH^k(M)$ for the subset of $H^k(M,P_k)$ of 
those cohomology classes represented by some $k$-admissible null-cobordant 
immersions. 
} 
\end{definiz} 
 
\begin{lem} $NEH^k(M)$ is a subgroup of $H^k(M,P_k)$. 
 
\end{lem} 
 
\begin{dimostraz} Given $\alpha$ and $\beta$  cohomology 
classes represented by $k$-admissible null-cobordant immersions it is 
obvious that $\alpha+ \beta$ is represented the same way. 
 
As for $-\alpha$, let $f$ be a $k$-admissible immersion such that 
$\chi^k_f=\alpha$, and let $g$ be a cobordism to the empty set of $f$. 
Let $g'$ be the cobordism between $f$ in $M\times\{1\}$ and the empty 
set in $M\times \{0\}$ obtained by composing $g$ with the reflection of 
$I$ given by $t\mapsto 1-t$. 
In a single $e^k\times I$ consider the following construction. Put a 
representative immersion of $-\alpha(e^k)$ in $e^k\times\{0\}$, put in 
$e^k\times \{1/3\}$ the same representative plus two copies of 
$f\cap{e^k}$ slightly isotoped, and in $e^k\times \{2/3\}$ a single copy of $f\cap{e^k}$. 
Then fill $e^k\times [1,1/3]$ with $-\alpha(e^k)\times[1,1/3]$ plus two copies 
of $g'\cap(e^k\times I)$ rescaled of 1/3, 
fill $e^k\times [1/3,2/3]$ with a cobordism between $-\alpha(e^k)$ plus 
a copy of $f\cap{e^k}$ and the empty set and with $f\cap{e^k}\times [1/3,2/3]$, and finally fill $e^k\times [2/3,1]$ 
with 
$g'\cap(e^k\times I)$ (rescaled of 1/3). Remark that this immersion 
is not a cobordism between $-\alpha(e^k)$ and $\alpha(e^k)$, since $g$ 
and $g'$ possibly intersect $(\partial e^k)\times I$. 
Now remark that the collection of immersions 
so defined glue together to a cobordism in $M_k\times I$, that 
restricted to $M_k\times \{0\}$ represents $-\alpha$. 
Consider on $M\times \{1\}$ the empty immersion and in $M_k\times I$ 
the collection of immersions defined before. This cobordism can be 
completed by lemma~\ref{cofibration} 
to  a cobordism $\tilde{g}$ between an immersion 
$\tilde{f}=\tilde{g}\cap(M\times\{0\})$ representing $-\alpha$ and the 
empty set, hence $-\alpha\in NEH^k(M)$.~$\square$ 
\end{dimostraz} 
 
The following lemma, that provides the technical step of the previous 
proof, will be repeatedly used in this section. It is an easy 
algebraic-topological argument, that has however an important geometric 
interpretation. 
 
\begin{lem}\label{cofibration} An immersion $g$ traced in $M\times \{0\}\cup M_s\times I$ 
extends to an immersion traced in the whole of $M\times I$. Moreover if 
$g$ is $l$-admissible, for $l\leq s$, the extension is still 
$l$-admissible. 
\end{lem} 
 
\begin{figure}[h] 
\begin{center} 
\includegraphics[width=8cm]{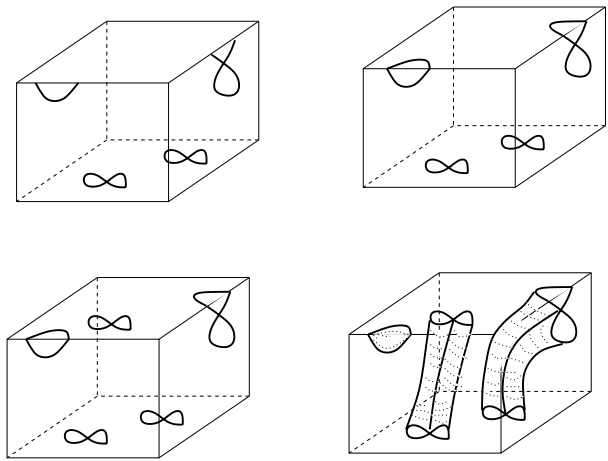} 
\caption{Extending an immersion to the interior of a cube with a free face}\label{completamento} 
\end{center} 
\end{figure}

\begin{dimostraz} If we see an immersion traced in $M\times \{0\}\cup M_s\times 
I$ as a continuous map from $M\times \{0\}\cup M_s\times I$  to $\QP$ 
this immediately follows from the fact that $(M,M_s)$ is  a 
cofibration. This construction has however an easy geometrical 
interpretation, see figure~\ref{completamento}. 
Consider a $(s+1)$-cube $e^{s+1}$. We want to extend the immersion traced 
on 
$g\cap 
(e^{s+1}\times\{0\}\cup \partial e^{s+1}\times I)$ to $e^{s+1}\times 
I$. First remark that the resulting immersion on the boundary of the 
cube at the top is null-cobordant, since it bounds in the disk $e^{s+1}\times\{0\}\cup \partial e^{s+1}\times 
I$, hence a null-cobordism can be traced on the cube at the top. The 
resulting immersion in $\partial(e^{s+1}\times I)$ represents an 
element of $P_{s+1}$. Up to adding (in the interior of the cube at the 
top) another immersion we can assume this element is trivial, hence 
  $g$ can be 
extended to the interior of $e^{s+1}\times I$. Recursively $g$ is extended to 
$M\times I$. 
 
The second statement follows obviously from the construction.~$\square$

\end{dimostraz} 
 
We prove that $\chi^k$ is a well-defined invariant in the group 
$H^k(M,P_k)/NEH^k(M)$.

\begin{prop}\label{welldefined} Let $f$ and $f'$ be $k$-admissible cobordant immersions. Then 
$\chi^k_f-\chi^k_{f'}$ represents an element in $NEH^k(M)$, hence 
 $\chi^k$ is a well-defined cobordism invariant in the group 
$H^k(M,P_k)/NEH^k(M)$. 
\end{prop} 
 
\begin{dimostraz} 
The 
technical step is  to modify $g$, far from $M\times\{0\}$, in such a way that it 
doesn't intersect $M_{k-2}\times I$. This is only possible, in general, up to 
modifying $g\cap M\times\{1\}$. 
By proposition~\ref{oneadmissible} we 
might assume that $g$ is 2-admissible. 
Assume by a recurrence hypothesis that $g$ is a $s$-admissible 
cobordism between $f$ and $f'+f_0$, where $f_0$ is a 
(null-cobordant) $s$-admissible immersion.

We want to deform $g$ in such a way that $g(e^{s-1})$ becomes trivial 
for any $e^{s}$, in order to apply lemma~\ref{sadmissible}. 
We build up an auxiliary immersion $\tilde{g}_s$ in $M\times I$. 
On the vertical cells $e^{s-1}\times I$ put a copy $g(e^{s-1})$, and on the 
bottom cells $e^{j}\times \{0\}$ fix the empty immersion for any $j$. 
Consider then 
the resulting collection of immersions as an immersion in 
$M\times\{0\}\cup M_{s-1}\times I$ and extend it to a cobordism 
$\tilde{g}_s$ by means of lemma~\ref{cofibration}. Then  $\tilde{g}_s$ satisfies the 
following properties: 
\[ 
\begin{array}{ll} 
\tilde{g}_s(e^{s-1})=g(e^{s-1}),&\\ 
\tilde{g}_s\cap M\times\{0\}=\emptyset,\\ 
\tilde{g}_s\cap M_{s-2}\times I=\emptyset. 
\end{array} 
\] 
Remark now that there exists an integer $a$ (the order of $P_{s}$ minus 1) such that $g+a\tilde{g}_s$ 
has the property that $(g+a\tilde{g}_s)(e^{s-1})$ is trivial for any $e^{s-1}$. Call again 
$g$ this cobordism and apply lemma~\ref{sadmissible}. The resulting 
cobordism $g'$ between $f$ and $g'\cap M\times\{1\}$ is then 
$(s+1)$-admissible, and $g'\cap M\times\{1\}$ is of the form 
$f'+f_0'$ where $f_0'$ is (null-cobordant and) $(s+1)$-admissible. 
 
Repeat this construction until $s=k-1$ and the resulting $k$-admissible cobordism, 
that we again call $g$, is such that 
\[ 
g\cap M\times\{1\}=f'+f_0, 
\] 
where $f_0$ is (null-cobordant and) $k$-admissible. 
Then by proposition~\ref{condizione} 
$\chi^k_{f'+f_0}$ differs from $\chi^k_f$ by a 
coboundary. But since 
$$\chi^k_{f'+f_0}=\chi^k_{f'}+\chi^k_{f_0},$$ 
one gets  the claim.~$\square$

\end{dimostraz}

The invariant 
\[ 
\chi^k:F^k\longrightarrow \frac{H^k(M,P_k)}{NEH^k(M)}, 
\] 
defined by these propositions will be called the {\em $k$-th cohomological 
invariant}.

\noindent Remark 
that if $f$ is actually in $F^{k+1}$ then $\chi^k_f$ is 
trivial, hence we are left with a well-defined homomorphism 
from $F^k/F^{k+1}$ to $H^k(M,P_k)/NEH^k(M)$. This homomorphism is in fact 
injective. 
 
\begin{prop}\label{nucleo} Let $M$ be an $n$-manifold. For any 
$k=1,\dots,n$ the kernel of the $k$-th cohomological invariant is 
$F^{k+1}$. 
\end{prop} 
 
\begin{dimostraz} Fix a cubulation in $M$. Let $f$ be a $k$-admissible immersion such that 
$\chi^k_f=0\in H^k(M,P_k)$. This means there is an element $\gamma\in 
C^{k-1}(M,P_k)$ such that $\delta\gamma=\chi^k_f$. One builds up 
a cobordism $g$ 
between $f$ and a $(k+1)$-admissible immersion. 
 
Consider the standard cubulation of $M\times I$ associated to the given 
cubulation of $M$. Put $f$ in the bottom $M\times \{0\}$. 
For any $(k-1)$-cell of $M$, say  $e^{k-1}$,  put in the 
vertical cell $e^{k-1}\times I$ 
(with the orientation that induces on $e^{k-1}$ 
the opposite of its orientation) the element $\gamma(e^{k-1})\in 
P_k$. One fixes also the cobordism on the top $M_k\times \{1\}$. 
Consider a $k$-cell $e^k$ of $M$. One defines the cobordism on $e^k\times 
I$ by remarking that (from the definition of $\gamma$) 
the union of all immersions already defined in 
$\partial (e^k\times I)$ is 
null-cobordant.  One can choose therefore a 
cobordism to the empty set. This can be done recursively on the whole 
of $\Mk\times I$. Remark that the resulting immersion $g$ does not 
intersect $\Mk\times \{1\}$.

Now $g$ is defined on $M\times\{0\}\cup M_k\times I$, and  applying 
lemma~\ref{cofibration} gives rise to a cobordism $g$ in $M\times I$, 
that provides a  cobordism  between $f$ and $g\cap(M\times 
\{1\})$; and this last does not intersect the $k$-skeleton, by 
construction. 
 
A similar construction can be performed if $\chi^k_f\in NEH^k(M)$. 
Let  $f_0$ be a null-cobordant map such that 
$\chi^k_{f_0}\sim\chi^k_f$. Put $f$ in the whole bottom $M\times \{0\}$, 
$f_0\cap M_k$ in the intermediate $M_k\times \{1/2\}$ 
and trace on the vertical $M_k\times [1/2,1]$ the intersection of a cobordism to the 
empty set of $f_0$. The $(k-1)$-cochain that cobounds $\chi^k_{f_0}$ and 
$\chi^k_f$ provides as before a cobordism between $f\cap M_k$ and 
$f_0\cap M_k$, which  we put in the vertical $M_k\times[0,1/2]$. 
Over all  this is a 
cobordism between $f\cap M_k$ and the empty set and so 
it  can be extended 
by lemma~\ref{cofibration} to a cobordism $g$ between $f$ and a map 
$g\cap(M\times\{1\})$.  By construction the last one 
does not intersect the $k$-skeleton.~$\square$ 
 
\end{dimostraz} 
 
\begin{coroll} Let $M$ be an $n$-manifold. For any 
$k=1,\dots,n$ the induced homomorphism 
\[ 
\tchi^k:\frac{F^k}{F^{k+1}}\longrightarrow \frac{H^k(M,P^k)}{NEH^k(M)}, 
\] 
is injective. 
\end{coroll}

This corollary shows that the power of these new invariants is 
considerable.  Indeed they 
describe the {\em graded group of $N(M)$ associated to the filtration} 
\[ 
gr(N(M))=F^1/F^2\times\dots \times F^{n-1}/F^n\times F^n, 
\] 
as a subgroup of $H^1(M,P_1)/NEH^1(M)\times\dots\times H^n(M,P_n)/NEH^n(M)$, that is: 
 
\begin{teo}\label{iniettivo} The cohomological invariants induce an injective homomorphism 
\[ 
\tchi:gr(N(M))\longrightarrow H^1(M,P_1)/NEH^1(M)\times\dots\times 
H^n(M,P_n)/NEH^n(M).~\square 
\] 
\end{teo}

\noindent We end this section with 
an important remark. 
The cohomological invariants reduce, under suitable 
hypothesis, to the restriction of James-Hopf invariants. These are 
classical cobordism invariants, see~\cite{JamOST} and 
\cite{KoSSIH}. 
 
\begin{definiz}{\rm Let $M$ be a $n$-manifold and $(F,f)$ a 
generic codimension-one  immersion. 
For $i=1,\dots,n$ consider the locus of $(n-i)$-tuple points 
of $f$, that is, the points of $M$ that have a number of preimages 
equal or bigger than $n-i$. This set is in fact a $i$-cycle modulo 2, whose homology class 
is invariant up to cobordism. We denote $JH_i(f)\in H_i(M,\zd)$ this 
class and call it {\em $i$-th James-Hopf invariant}. 
} 
\end{definiz} 
 
\noindent These invariants are particularly meaningful for those $k$ such that 
$P_k$ is non-trivial. 
For example given a codimension-$k$  embedded submanifold 
$S$ of an $n$-manifold $M$ 
such that its normal bundle is reducible to the symmetry group of an 
element $f\in P_k$ 
satisfying 
$\theta_k(f)\neq0$ then 
  the following holds 
\[ 
JH_i(S\ltimes f)=\left\{ 
\begin{array}{ll} 
[S]\in H_{n-k}(M,\zd)&\text{ if } i=n-k\\ 
0&\text{ otherwise.} 
\end{array} 
\right. 
\]

 
\begin{prop} Assume that $k$ is such that 
$\theta_k:P_k\longrightarrow\zd$ is the reduction modulo 2 of 
$P_k$. Then for any $n$-manifold the $(n-k)$-th 
James-Hopf invariant restricted to $F^k$ is the Poincar\'e dual to the reduction modulo 2 of 
$\chi^k$. 
\end{prop} 
 
\begin{dimostraz} Let $f$ be a $k$-admissible immersion generic and 
  transverse to the decomposition of $M$. We denote by $PD$ the 
  Poincar\'e duality isomorphism. Then for any 
$k$-cell $e^k$ of $M$ the number 
$PD(JH_{n-k}(e^k))$ is the number of $k$-tuple points of $f\cap e^k$, 
modulo 2, hence, by the hypothesis on $k$, is the reduction modulo 2 of 
$f\cap e^k$ as an element of $P_k$, that is $\theta_k(\chi^k_f(e^k))$. The following diagram then commutes 
\[ 
\xymatrix{ 
F^k\ar[r]^{\chi_k}\ar[d]_{JH_{n-k}|_{F^k}}&H^k(M,P_k)\ar[ddl]^{\theta_{k}^*}\\ 
H_{n-k}(M,\zd)\ar[d]_{PD}\\ 
H^k(M,\Z/2\Z) 
} 
\] 
and since $\theta_k^*$ is reduction modulo 2 in cohomology the claim follows.~$\square$ 
\end{dimostraz} 
 
\noindent This proves at once the following: 
\begin{prop} 
If $k$ is such that $\theta_k$ is an isomorphism then 
$NEH^k(M)=0$.~$\square$ 
\end{prop}

\section{Obstruction theory}\label{ostruzione} 
 
A more detailed study of the groups $NEH^k(M)$ is in order. From 
proposition~\ref{oneadmissible} 
and lemma \ref{sadmissible} one might guess that the vanishing of 
$NEH^k(M)$ is correlated with constructions that make the 
intersections $g(e)$ of a cobordism $g$ with the vertical walls of $M\times 
I$ null-cobordant. 
This was made possible for example in the proof of 
proposition~\ref{welldefined} by means of the construction of an 
auxiliary immersion $\tilde{g}$ with prescribed image on $M\times 
\{0\}\cup M_{s-1}\times I$.  As we saw the possibility of obtaining such an 
auxiliary immersion amounts, at an accurate 
analysis, to 
the fact that $(M,M_{s-1})$ is a cofibration. However in order to leave  the 
image of $g$ fixed also in $M\times \{1\}$ we need a more delicate 
construction, that will be developed in section~\ref{ostruzione} and 
applied in section~\ref{NEH}. 
 
A second obvious motivation for developing this 
theory is the computation of the image of $\chi^k$, that is, the 
subgroup of cohomology classes that are represented by an immersion.

In this section we describe the general obstructions for a cochain in 
$C^k(M,P_k)$ to 
be realizable as an immersion. The basic idea is a recursive 
construction. Given $\xi\in C^k(M,P_k)$ we first put in the interior of every 
$k$-cube an immersion representing $\xi(e^k)\in P_k$, then  try to extend 
this codimension-one immersion in $M_k$ to a codimension-one immersion 
in $M_{k+1}$.  This will be called the {\em first extendibility}. 
If this construction can be repeated until the $n$-th 
skeleton i.e. the immersion can be further extended to a 
{\em second extension}, and  so on 
then the original cochain is said to be {\em realizable} or {\em extendible}. 
If one can reach the $s$-th stage one says 
that the cochain is {\em $s$-extendible}. 
 
We adapt to this context Eilenberg's obstruction theory, see 
\cite{WhiEHT}, \S V.5. A cochain in $C^k(M,P_k)$ can be thought of as a 
map $\varphi_k$ defined from $M_k$ to $\QP$, that restricted 
to any $k$-cell is geometrically represented by an element of $P_k$. 
The problem to which we apply Eilenberg's theory is that of extending 
this map over the next skeleton. 
 
\subsection{A review of obstruction theory}

Given a $s$-simple space $Y$, a CW-complex $X$ and a map 
$\varphi:X_s\longrightarrow Y$ 
the {\em obstruction to extending $f$ to the $(s+1)$-skeleton} is a 
cochain $c^{s+1}(\varphi)\in C^{s+1}(X,\pi_s(Y))$, assigning to each 
$(s+1)$-cell $e^{s+1}$ the map $\varphi|_{\partial e^{s+1}}$. 
Its fundamental 
properties are stated in the following theorem (see \cite{WhiEHT}, \S V.5): 
 
\begin{teo}\label{properties} 
 
\begin{enumerate} 
\item $\varphi$ is extendible to $X_{s+1}$ if and only if  $c^{s+1}(\varphi)$ 
is the trivial cochain. 
\item $c^{s+1}(\varphi)$ is a cocycle, hence represents an element of 
  $H^{s+1}(X,\pi_s(Y))$. 
\item $\varphi|_{X_{s-1}}$ is extendible to $X_{s+1}$ if and only if 
$c^{s+1}(\varphi)$ is trivial in $H^{s+1}(X,\pi_s(Y))$. 
\end{enumerate} 
\end{teo} 
 
\noindent The problem of further extending $\varphi$ is 
codified in a sequence of 
obstruction maps. However for any extension there exists an obstruction 
cocycle, hence the obstruction to further extend $\varphi|_{X_{s-1}}$ 
becomes a set of cohomology classes. Assume that it is extendible to 
the $(s+l-1)$-th skeleton, and 
let $\o^{s+l}(\varphi)$ be the set of 
cohomology classes given by 
\[ 
\o^{s+l}(\varphi):=\{c^{s+l}(\varphi_{s+l-1})| 
\varphi_{s+l-1}\mbox{ is an extension of } \varphi|_{X_{s-1}} 
\mbox{ to the }(s+l-1)\mbox{-th skeleton}\} 
\subset H^{s+l}(X,\pi_{s+l-1}(Y)). 
\] 
\begin{teo}\label{extendability} Assume the map $\varphi|_{X_{s-1}}$ is extendible to the 
$(s+l-1)$-skeleton.  Then it is extendible to the 
$(s+l)$-skeleton if and only if the 0 class belongs to 
$\o^{s+l}(\varphi)$. 
 
\end{teo}

\subsection{Geometrical interpretation} 
 
We first recall that 
\[ 
[S^s,\QP]=P_s, 
\] 
hence the theory applies directly with coefficients in the groups 
$P_s$. 
 
Fixed a $k$-cochain $\xi\in C^k(M,P_k)$ we consider it as a map 
defined on $M_k$ taking values on 
$\QP$. 
Consider then the obstruction 
\[ 
\begin{array}{rccc} 
c^{k+1}(\xi):&C_{k+1}(M)&\longrightarrow& \pi_{k}(\QP),\\ 
&e^{k+1}&\mapsto&\xi(\partial e^{k+1}), 
\end{array} 
\] 
and remark that, since 
$\pi_{k}(\QP)=P_k$, the 
obstruction 
$c^{k+1}$ is nothing but the ordinary coboundary of cochains with 
coefficients in $P_k$.  Hence by property 1 of theorem~\ref{properties} 
it follows that: 
 
\begin{prop}$\xi$ is 1-extendible if and only if it is a 
cocycle. 
\end{prop} 
 
\begin{oss}{\rm This condition has the  geometrical interpretation 
that was already illustrated in figure~\ref{closed}. 
 } 
\end{oss} 
 
\noindent Define now $\o_k^{k+l}$, for $l\geq2$, 
to be the map that associates to $\xi\in H^k(M,P_k)$ 
the set $\o^{k+l}(\varphi_{\xi})\subset H^{k+l}(X,P_{k+l-1})$, $\varphi_{\xi}$ being any extension of $\xi$ 
to the $(k+1)$-skeleton, and define $\ker\o_k^s$ to be the subset of 
$H^k(M,P_k)$ of cocycles $\xi$ such that $\o^s(\varphi_{\xi})$ contains 
the trivial element of $H^{s}(X,P_{s-1})$. 
Remark that $\ker\o_k^{s-1}\subseteq \ker\o_k^{s}$, and that 
$\o_k^{s}$ is in 
fact only defined on $\ker\o_k^{s-1}$. 
From 
theorem~\ref{extendability} one obtains the following proposition: 
 
\begin{prop}\label{realizability}$\xi\in H^k(M,P_k)$ is $l$-extendible if and only if it 
belongs to  $\ker\o_k^{k+l}$. In particular 
it is realizable if and only if it belongs to $\ker\o_k^n$.~$\square$ 
\end{prop} 
 
\noindent Denote by $EH^k(M)$ the subgroup of $H^k(M,P_k)$ of 
extendible cocycles. 
Proposition~\ref{realizability} translates into: 
\[ 
EH^k(M)=\ker\o_k^n. 
\]

 
\begin{prop} Let $\xi\in EH^k(M)$; then the immersion $f_{\xi}$ 
realizing $\xi$ is a well-defined element in $F^k/F^{k+1}$. 
\end{prop} 
 
\begin{dimostraz} This follows from 
proposition~\ref{nucleo}. If $f'$ and $f$ both realize $\xi$ in 
particular they have the same $k$-th cohomological invariant $[\xi]\in 
EH^k(M)/NEH^k(M)$, so they differ by an element of 
$F^{k+1}$.~$\square$ 
 \end{dimostraz}

\section{Explicit computations} 
 
We apply the obstruction theory to the computation of both 
$EH^k(M)$ and 
$NEH^k(M)$ (the subgroup of $EH^k(M)$ of {\em null-extendible 
cocycles}  cocycles realizable as a null-cobordant $k$-admissible 
immersion). The computations prove theorem~\ref{main}. 
 
\subsection{Extendible cocycles}\label{EH}

The results are summarized in the following 
table
%
%
, where 
$N^k:=H^k(M,P_k)/EH^k(M)$ denotes the {\em non-extendible cocycles}.

\vspace{0.5cm} 
 
{\small 
 
\begin{center} 
 
\renewcommand{\arraystretch}{1.2} 
 
\begin{tabular}{|c|c|c|c|c|c|c|c|c|} \hline 
dim $M$& orientability &  $N^1$ & $N^2$ & $N^3$ & $N^4$ & $N^5$ & $N^6$ & $N^7$ \\ \hline 
2      & both & $0$&  $0$ & $\cdot$& $\cdot$ &$ \cdot$&$\cdot$&$\cdot$\\ \hline 
3      & both & $0$&  $0$ & $0$& $\cdot$ &$ \cdot$&$\cdot$&$\cdot$\\ \hline 
4      & orientable & $0$&  $H^2(M,\zd)/Q^2(M,\zd)$ & $0$& $0$ &$ \cdot$&$\cdot$&$\cdot$\\ \hline 
4      & non-orientable & $0$&  $0$ & $0$& $0$ &$ \cdot$&$\cdot$&$\cdot$\\ \hline 
5      & both & $0$&  $H^2(M,\zd)/Q^2(M,\zd)^{(*)}$ & $0$& $0$ &$ 0$&$\cdot$&$\cdot$\\ \hline 
6      & both & $0$&  $H^2(M,\zd)/Q^2(M,\zd)^{(*)}$ & $0$& $0$ &$ 0$&$0$&$\cdot$\\ \hline 
7      & both & $0$&  $H^2(M,\zd)/Q^2(M,\zd)^{(*)}$ & $0$& $0$ &$ 0$&$0$&$0$\\ \hline 
 
\multicolumn{9}{l}{ (*) Under the condition $M$ orientable and ${\rm Ext}(H_3(M,\Z),\zo)=0$.} 
\end{tabular} 
 
\renewcommand{\arraystretch}{1} 
\end{center}

} 
 
\vspace{0.5cm} 
 
\noindent In this table  $Q^2(M,\zd)$ denotes the subgroup of $H^2(M,\zd)$ defined by 
\[ 
Q^2(M,\zd)=\{x\in H^2(M,\zd); \, x\cup x=0\}. 
\] 
Remark that, if $M$ is a manifold, this is  the quadric of 
$H^2(M,\zd)$ associated with the intersection form. 
 
\noindent Moreover for any dimension $n$ one has: 
 
\vspace{0.5cm} 
\begin{center} 
{
\begin{tabular}{|c|c|c|c|c|c|}\hline 
dim $M$& orientability &  $N^1$ & $N^s$ & $N^{n-1}$ & $N^n$ \\ \hline 
$n$      & both & $0$&  $?$ & $0$& $0$ \\ \hline 
\end{tabular} 
} 
\end{center} 
\vspace{0.5cm}

\noindent We first prove the results from the last table. 
 
\begin{prop}\label{anyn} For any dimension $n$ and  for any 
$n$-manifold $M$ 
\[ 
EH^1(M)=H^1(M,P_1)=H^1(M,\zd), \] 
\[ EH^{n-1}(M)=H^{n-1}(M,P_{n-1}), \] 
\[ EH^n(M)=H^n(M,P_n)=\left\{\begin{array}{ll} 
P_n, &\text{ if $M$ is orientable,} \\ 
P_n/2P_n, &\text{ if $M$ is non-orientable.} \\ 
\end{array}\right.\] 
\end{prop} 
 
\begin{dimostraz} Since $P_1=\zd$ one can use Poincar\'e duality 
with coefficients $P_1$ in both orientable and non-orientable 
context. Now represent the  Poincar\'e dual to $\xi$ 
by an embedding  and remark that this embedding 
realizes $\xi$. 
That $n$ and $(n-1)$-cohomology classes extend follows immediately from 
obstruction theory. ~$\square$ 
 

\end{dimostraz}

\begin{oss}{\rm This result can be interpreted geometrically. 
That every $n$-class is extendible follows easily from the fact that taking a 
representative cocycle $\xi$ a putting in any $n$-cube $e^n$ an immersion representing 
$\xi(e^n)$ already realizes the cocycle. 
 
Let now $M$ be an orientable $n$-manifold, let $\xi$ be a cohomology 
class in $H^{n-1}(M,P_{n-1})$, and consider its Poincar\'e dual $PD(\xi)$. 
By the universal coefficient theorem $PD(\xi)$ can be 
thought of as an element of $H_1(M,\Z)\otimes P_{n-1}$, hence as a 
combination of the type $\sum \gamma_i\otimes f_i$ with $\gamma_i\in 
H_1(M,\Z)$ and $f_i\in P_{n-1}$. 
To each 
$\gamma_i\otimes f_i$ we associate the following immersion. Take a simple 
closed loop representing $\gamma_i$. By decorating $\gamma_i$ with 
$f_i$, which is always possible since the normal bundle of $\gamma_i$ 
is trivial, 
 one obtains an immersion realizing the Poincar\'e dual 
of $\gamma_i\otimes f_i$. Obviously the sum of such immersions realizes 
$\xi$. 
 
Let $n$ be such that $P_{n-1}$ is either trivial or $\zd$. 
Then if $M$ is a non-orientable $n$-manifold the same construction 
applies. Indeed, 
these conditions on $P_{n-1}$ both mean that any 
immersion in $\R^{n-1}$ (up to a cobordism) admits  a reflection in its symmetry group, since $f=-f$. 
Hence 
this construction applies, since also curves with 
non-orientable normal bundle can be always decorated with immersions in 
$P_{n-1}$.

} 
\end{oss}

We are ready now to prove the results from the main table. 
Let us first concentrate on codimension-two cocycles. If $k=n-2$ there 
is 
only one obstruction, namely $\o_{n-2}^n$. This lives 
in the cyclic group $H^n(M,P_n)$. One  describes now 
$\o_{n-2}^n(\xi)$ explicitly as an element of $P_n$ or of 
$P_n/2P_n$, depending on 
$M$ being orientable or not.

\begin{lem}\label{almostextendable} Let $M$ be any $n$-manifold. Every $\xi\in H^{n-2}(M,P_{n-2})$ is extendible to an 
immersion defined in $M\setminus int(e^n)$, where $e^n$ is an $n$-ball. 
 
\end{lem} 
 
\begin{dimostraz} Perform a first extension $f$ of $\xi$ to $M_{n-1}$. 
Remark that $M\setminus int(e^n)$ collapses 
simplicially on a subset $S_{n-1}$ of $M_{n-1}$. Fix a way of building 
up $M\setminus int(e^n)$ from this subset, that is, order the set of 
$n$-cells in such a way that $e^n_{(1)}$ is attached to $S_{n-1}$ and 
$e^n_{(i)}$ 
is attached to $S_{n-1}\cup \bigcup_{j=1}^{i-1} e^n_{(j)}$. For any $i$ 
call {\em free face} of $e^n_{(i)}$ the one to which first a cell of $M$ will be 
attached. Remark that every $n$-cell but $e^n$ has a free face. 
When 
attaching the first $n$-cell $e^n_{(1)}$ extend the immersion 
$f\cap{\partial e^n_{(1)}}$ 
 this way. If $f\cap{\partial e^n_{(1)}}$ is trivial in $P_{n-1}$ then the 
extension is the cobordism to the empty set, if it is non-trivial, see 
figure~\ref{nonextendable}, 
\begin{figure}\label{nonextendable} 
\hspace{7cm}\psfig{figure=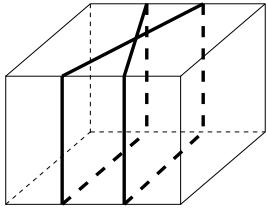,width=3.5cm} 
\caption{An immersion not extendible over the ball} 
\end{figure} 
 then add a representative of $-f\cap{\partial 
e^n_{(1)}}$ on the free face of $e^n_{(1)}$. Call again $f$ the new 
extension of $\xi$, and perform recursively the same construction. At 
the end one is left with $f$ defined on $M\setminus int(e^n)$.~$\square$ 
 
\end{dimostraz} 
 
\begin{lem}\label{independence} Let 
$M$ be an  $n$-manifold, $\xi\in H^{n-2}(M,P_{n-2})$ and 
$f$ be any extension of $\xi$ to $M\setminus int(e^n)$. 
If $M$ is orientable then $f\cap{\partial e^n}\in P_{n-1}$ depends only on 
$\xi$.  If $M$ is non-orientable  then 
$f\cap{\partial e^n}\in P_{n-1}/2P_{n-1}$ depends only on 
$\xi$. 
\end{lem} 
 
\begin{dimostraz} Given any extension $f$ of $\xi$ to $M_{n-1}$, 
the set of extensions modulo cobordism relative to $M_{n-2}$ is acted on 
by $C_{n-1}(M,P_{n-1})$, by the action $(\alpha * f)\cap{e^n}=\alpha(e^n)\cup 
(f\cap{e^n})$. This action is transitive. Remark that if two extensions $f$ 
and $f'$, both extends to $M\setminus int(e^n)$, then their difference 
$\alpha$ must be such that for any $n$-cell $e'\neq e$ it holds 
$\alpha(\partial e')=0$. 
 
Assume now that $M$ is orientable. Then at the cochain level $\partial 
e=\sum_{e'\neq e} \partial e'$, so $\alpha (\partial e)=0$, and hence 
$f\cap{\partial e}=f'\cap{\partial e}$. 
 
If $M$ is non-orientable the equation $\partial 
e=\sum_{e'\neq e} \partial e'$ only holds modulo 2, hence one can say 
that $\alpha (\partial e)\in 2P_{n-1}$.

Assume now that $\xi$ is represented by a different cocycle, hence by a 
different immersion in $M_{n-2}$. Since the two cocycles are 
cohomologous the two immersions are cobordant, hence there exists a 
cobordism in $M_{n-2}\times I$ between the two representatives. This 
cobordism can be extended to a cobordism between extensions to 
$M_{n-1}$, since any cube $e^{n-1}\times I$ has a free face, say, $e^{n-1}\times 
\{1\}$, 
 and in an analogous way to a cobordism between two extensions $f$ and 
$f'$ to $M\setminus int(e^n)$.  This proves that $f\cap{\partial e^n}$ 
and $f'\cap{\partial e^n}$ are cobordant. 
 
That $f\cap{\partial e}$ does not depend neither on $e$ nor on the process of 
collapsing is then straightforward, hence the claim.~$\square$ 
\end{dimostraz} 
 
\begin{prop}\label{codimensiontwo} Let $M$ be an  $n$-manifold and $\xi\in 
H^{n-2}(M,P_{n-2})$. Then $\xi\in EH^{n-2}(M)$ if and only if, given any extension $f$  of $\xi$ to $M\setminus int(e^n)$, 
\[ 
f\cap{\partial e}=0 \in\left\{ 
\begin{array}{ll} 
P_{n-1}, & \text{\ if } $M$ \text{ is orientable,}\\ 
P_{n-1}/2P_{n-1},& \text{\ if }$M$ \text{ is non-orientable}. 
\end{array} 
\right. 
\] 
\end{prop} 
 
\begin{dimostraz} The case $M$ orientable follows immediately from the preceding lemmas. As for the 
non-orientable case recall the proof of proposition~\ref{Q3} and proposition~\ref{pari} and remark that if $f\cap{\partial 
e}\in 2P_{n-1}$ then its opposite $g$ bounds in $M\setminus int (e)$ 
a cobordism not intersecting $M_{n-2}$.~$\square$ 
 
\end{dimostraz} 
 
\begin{coroll} Let $n$ be such that $Q_{n-1}=0$. Then for any 
$n$-manifold 
\[ 
EH^{n-2}(M)=H^{n-2}(M,P_{n-2}). 
\] 
More generally, if $n$ is such that $Q_{n-1}\subseteq 2P_{n-1}$, 
then the same holds true for 
non-orientable $n$-manifolds. 
\end{coroll} 
 
\begin{dimostraz} This immediately follows from proposition~\ref{codimensiontwo} and the fact that for any 
extension $f$ of $\xi$ to $M\setminus int (e^n)$ the immersion 
$f\cap{\partial e}$ belongs to $Q_{n-1}$.~$\square$

\end{dimostraz}

\noindent This yields the claimed  values for $N^{n-2}$ in all cases 
 but for orientable $4$-manifolds. In 
this case a geometric construction is in order. 
 
\begin{prop}\label{twofour} Let $M$ be an orientable 4-manifold and $\xi\in 
H^2(M,\zd)$. Then $\xi\in EH^2(M)$ if and only if $\xi\cup \xi=0$. 
\end{prop}

\begin{dimostraz} Assume first that $\xi\cup\xi=0$, that is, $\xi\in 
Q^2(M,\zd)$. Take a smoothly embedded representative $F$ of $PD(\xi)$, and 
take a generic normal field $\nu$ to $F$ in $M$. The hypothesis on $\xi$ 
implies that $\nu$ has an even number of isolated, simple zeroes $z_1,\dots,z_{2s}$. Around each zero $z_i$ take a 
small disk $D_i$ in $F$ such that $\nu|_{\partial D_i}$ has degree 1 or -1. 
Then cut off all the disks and connect the remaining holes in pairs 
with tubes which are contained in a tubular neighborhood of $F$. The 
resulting surface $F'$ still represents $PD(\xi)$ and admits a nowhere 
vanishing normal 
field of directions. Since the group of symmetries of the 8 in $\R^2$ 
is equal to 
the group of symmetries of a line, the existence of the field of 
directions means that it is possible to decorate $F'$ with 8's. The 
resulting codimension-one immersion is in $F^2$ and has second 
cohomological invariant equal to $\xi$. 
 
On the opposite direction, assume by absurd that there exists 
$\xi\in EH^2(M)$ with $\xi\cup\xi\neq0$. 
Consider  an embedded surface representing $PD(\xi)$, then 
choose a normal field to 
$F$. Up to changing $F$ in its homology class we can assume as in 
the previous step that $F$ admits a normal field of directions with 
a single isolated degenerate point $z$. 
Let $e^4$ be a 4-disk around $z$ in $M$.  One can 
extend $\xi$ to an immersion defined in $M\setminus int (e^4)$ by 
decorating $F$ with 8's following the normal field, call $f$ this 
immersion. 
Now $f\cap{\partial e^4}$ is a non-trivial element of 
$Q_3=2P_3$.  If it was  trivial, the normal field 
of directions could be extended to the whole of $F$, 
which is not possible. By a local 
analysis, we can reduce ourselves to the situation of 
remark~\ref{bounding_in_orientable}, hence 
one obtains that $f\cap{\partial e^4}$ 
is in fact the element $4\in P_3$. 
But 
from proposition~\ref{codimensiontwo}, since $M$ is orientable,  $\xi\notin 
EH^2(M)$.~$\square$

\end{dimostraz}

\begin{oss}{\rm We showed in proposition~\ref{Q3} that $Q_3=2P_3$ and in remark~\ref{bounding_in_orientable} that 
$4P_3$ contains immersions that bound in an orientable manifold. 
 The 
theory of this section shows that $4P_3$ is the subgroup of immersions 
bounding in an orientable manifold, that is, immersions 
with invariants 
2 and 6 do {\em not} bound in any orientable manifold. 
} 
\end{oss} 
 
This settles the table for $n\leq 4$. 
One  proves now that all obstructions involved in the table are 
trivial except for $\o_2^4$. 
This fact is due to the particular  properties of 
$P_s$ and $Q_s$ for $s\leq6$ and 
$s\neq4$. 
In general, the more the 
groups $P_s$ are simple the more the extensions $\o_k^{s+1}$ are easy 
to compute. The easiest case is of course $s=5$. 
 
\begin{prop} For $k\leq 4$ the extension $\o_k^6$ is 
trivial.~$\square$ 
\end{prop} 
 
\noindent The easiest next step is  a property of the first obstruction 
$\o_k^{k+2}$ for some values of $k$. 
 
\begin{prop}\label{theta} Let $k$ be such that $\theta_{k+1}$ is an 
isomorphism. Then the first 
obstruction $\o^{k+2}_{k}$ is trivial. 
 
\end{prop} 
 
\begin{dimostraz} For $\xi\in H^k(M,P_k)$ consider any first extension 
$f$ to $M_{k+1}$. Remark that $f\cap{e^{k+1}}$ cannot be considered 
as an 
element of $P_{k+1}$, since $f\cap \partial e^{k+1}$ is not trivial. 
However the number of $(k+1)$-tuple points modulo 2 of $f\cap{e^{k+1}}$  is 
well-defined. 
Then let $\alpha\in C^{k+1}(M,P_{k+1})$ be the cochain that associates 
to $e^{k+1}$ this number. Remark that the composition with $\theta_{k+1}$ 
induces a natural isomorphism between $C^{k+1}(M,P_{k+1})$ and 
$C^{k+1}(M,\zd)$. So we can consider $\alpha * f$, with the action 
defined in the proof of lemma~\ref{almostextendable}. This immersion 
extends $\xi$ and is extendible to $M_{k+2}$, since for any 
$(k+2)$-cell $e^{k+2}$ 
\[ 
\theta_{k+1}((\alpha * f)\cap{\partial e^{k+2}})= \theta_{k+1}(f\cap{\partial e^{k+1}})+\sum_{e\in\partial 
e^{k+1}}\alpha(e)=0, 
\] 
that is, $(\alpha * f)\cap{\partial e^{k+2}}=0\in P_{k+1}$.~$\square$ 
 
\end{dimostraz} 
 
The triviality of 
almost all of the obstructions involved  in the table follow then from generalizing the 
previous results.

The proof of proposition~\ref{theta} actually extends to the following 
result, that is in fact the strongest triviality result 
in this section. Given an 
immersion however traced on a $s$-skeleton, if $\theta_s$ is an 
isomorphism then in each $s$-cube one can force the parity of $s$-tuple 
points to be even, and since any $(s+1)$-cube has an even number of 
faces, 
this permits extension to the $(s+1)$-skeleton. 
 
\begin{teo}\label{forte} Let $s$ be such that $\theta_s$ is 
an isomorphism. Then the 
obstruction $\o^{s+1}_{k}$ is trivial, for any $k\leq s-1$. 
 
\end{teo} 
 
\begin{dimostraz} The action defined in the proof of 
proposition~\ref{theta} can be defined on the set of extensions from 
any skeleton to the following one, hence the proof applies.~$\square$ 
 
\end{dimostraz}

\noindent On the other side, the construction of the unique 
obstruction for codimension-two cocycles can be performed 
in a more general context. Specifically, by  an easy adaptation of the 
arguments of lemmas~\ref{almostextendable} 
and \ref{independence} one proves the following proposition: 
 
\begin{prop} Let $n$ be such that $Q_{n-1}=0$. Then for any 
$n$-manifold and any $k\leq n-2$ the last obstruction 
$\o_k^n$ is trivial. 
If $n$ is such that $Q_{n-1}\subseteq 2P_{n-1}$ then the same holds for any 
non-orientable $n$-manifold.~$\square$ 
\end{prop} 
 
We are then left to study $N^2$ for $n$-manifolds with $n>4$, since 
$\o^4_2$ is the only nontrivial obstruction involved in the table. 
 
\begin{teo}\label{fiveseven} Let $M$ be a closed orientable $n$-manifold such that 
\[ 
{\rm Ext}(H_3(M,\Z),\zo)=0. 
\] 
Then 
\[ 
EH^2(M)\subseteq Q^2(M,\zd)\subseteq \ker\o_2^4. 
\]

\end{teo}

\begin{dimostraz} 
One shows first that  ${\rm Ext}(H_3(M,\Z),\zd)=0$ implies  $EH^2(M)\subseteq 
Q^2(M,\zd)$. 
This condition on ${\rm Ext}$ is in fact a consequence of ${\rm Ext}(H_3(M,\Z),\zo)=0$. 
Let $\xi\in EH^2(M)$.  The hypothesis on $M$ implies by 
the universal coefficient theorem that $H^4(M,\zd)={\rm Hom}(H_4(M,\Z),\zd)$. 
We want to prove that $\xi\cup\xi=0$ by showing that for all $c\in 
H_4(M,\Z)$ one obtains $(\xi\cup\xi)(c)=0$. Fix a 
class $c\in H_4(M,\Z)$, and represent it by an 
embedded orientable 4-submanifold $F$ (see \cite{ThoQPG}, theorem 
II.27). 
One  obtains easily $i^*(\xi)\cup 
i^*(\xi)=0$, where $i$ is inclusion of $F$ in $M$. 
In fact  $i^*(\xi)$ is 
an extendible element of $H^2(F,\zd)$, and the claim follows from 
the 
characterization of proposition~\ref{twofour}. Hence $0=i^*(\xi)\cup 
i^*(\xi)=i^*(\xi\cup\xi)=PD(\xi\cup\xi)(c)$. 
 
Now we show that ${\rm Ext}(H_3(M,\Z),\zo)=0$ implies that 
$Q^2(M,\zd)\subseteq \ker\o_2^4$. Let $\xi\in Q^2(M,\zd)$, then 
$\o^4_2(\xi)\in H^4(M,\zo)={\rm Hom}(H_4(M,\Z),\zo)$. For any $c\in 
H_4(M,\Z)$ take as before a representative $F$ orientable and embedded in $M$. 
Remark that $i^*(\xi)\cup i^*(\xi)=0$, hence $\o_2^4(i^*(\xi))=0\subset 
H^4(F,\zo)$. By functorial properties of obstruction cocycles 
(see \cite{WhiEHT} page 230) $\o_2^4(i^*(\xi))=i^*(\o_2^4(\xi))$, 
and this last is the set $\o_2^4(\xi)\subset 
H^4(M,\zo)={\rm Hom}(H_4(M,\Z),\zo)$ evaluated on $c\in H_4(M,\Z)$ and 
composed with Poincar\'e duality. 
We proved therefore that  this set contains $0$ 
for any $c\in H_4(M,\Z)$, hence $\o_2^4(\xi)$ contains 
the trivial cocycle and so  $\xi$ extends over the 4-skeleton.~$\square$ 
 \end{dimostraz}

 
%

\begin{coroll} If $M$ is a closed  orientable $n$-manifold, $n\leq 7$, 
and ${\rm Ext}(H_3(M,\Z),\zo)=0$ then 
\[ 
EH^2(M)=Q^2(M,\zd). 
\] 
\end{coroll} 
 
\begin{dimostraz} This is because for $n\leq7$ the obstruction 
  $\o_2^4$ is the only one that 
has not been proven to be trivial yet, hence $EH^2(M)=\ker\o_2^4$. So from theorem~\ref{fiveseven} 
we obtain the claim.~$\square$ 
 
\end{dimostraz}

\noindent The results summarized in the table are then proved. 
 
\noindent Finally remark that the crucial property of $\theta_s$ 
that makes 
proposition~\ref{theta} and theorem~\ref{forte} work 
can be abstracted to the following  definition: 
 
\begin{definiz}{\rm For any $f\in 
P_{s-1}$ denote by $N(B^s,f)$ the group of immersions $g$ in $B^s$ such 
that $g\cap\partial B^s$ are cobordant to $f$. 
 We say $s$ is {\em simple} if for any $f\in 
P_{s-1}$ 
there exists an 
isomorphism $i_f:N(B^s,f)\rightarrow P_s$ such that 
\[ 
\xymatrix{ 
N(B^s,f)\times N(B^s,-f)\ar[r]\ar[d]_{i_f\times i_{-f}}&N(S^s)\ar[d]^{\sim}\\ 
P_s\times P_s\ar[r]_{\cdot}&P_s. 
} 
\] 
} 
\end{definiz} 
 
\noindent Hence the following theorem holds: 
 
\begin{teo} For simple $s$ the 
obstruction $\o_k^{s+1}$ is trivial for any $k\leq s-1$.~$\square$ 
 
\end{teo}

\subsection{Null-extendible cocycles}\label{NEH}

We apply obstruction theory to computations concerning $NEH^k(M)\subseteq 
EH^k(M)$, the subgroup of cocycles that are extendible to 
null-cobordant immersions. We actually obtain the following satisfactory 
statement: 
 
\begin{teo} Let $M$ be an $n$-manifold, with $n\leq7$. Then for any 
$k=1,\dots,n$ 
\[ 
NEH^k(M)=0. 
\] 
 
\end{teo} 
 
\begin{dimostraz} We apply proposition~\ref{condizione}. We have to prove 
that, for any $k\leq n\leq7$, given a cobordism $g$ between 
$k$-admissible immersions,  one can obtain from $g$ a $k$-admissible 
cobordism. By lemma~\ref{sadmissible} it is enough to show recursively 
that if $g$ is $s$-admissible, for $s\leq k-1$, then one can obtain from 
it a $s$-admissible cobordism 
$g'$ such that for any $(s-1)$-cube $e^{s-1}$ of $M$ the 
intersection $g'(e^{s-1})$ is a trivial 
element of $P_{s}$.  Recall that $g'(e^{s-1})$ actually represents an 
element of $P_{s}$ since $g'$ is $s$-admissible. We already know by 
proposition~\ref{oneadmissible} that any cobordism $g$ between 2-admissible 
immersions 
can be considered to be 2-admissible. 
 
Assume then that $f$ and $f'$ are 3-admissible immersions, and let $g$ be a 
2-admissible cobordism between them. For any 1-edge $e^1$ of $M$ the 
immersion $g(e^1)$ represents a well-defined element of $P_2=\zd$. We 
define a 2-cocycle in $M\times I$ this way 
\[ 
\Lambda(e)=g\cap e. 
\] 
This means that if $e$ is horizontal $\Lambda(e)=0$ and if $e$ is 
vertical of the form $e^1\times I$ then $\Lambda(e)=g(e^1)$. 
Remark  that $\Lambda$ is closed as a cocycle in $M\times(0,1)$, since 
$f$ and $f'$ 
being 3-admissible implies that 
$g\cap(M_2\times \{0,1\})$ is empty. 
 
We claim that $\Lambda$ is an extendible cocycle (with compact support) in 
$M\times(0,1)$.  If it is so, then any associated 
immersion $h$ in $M\times(0,1)$ is such 
that $g'=g\cup h$ is a 2-admissible cobordism such that $g'(e^1)=0$ 
for all $e^1$, hence can be deformed to a 3-admissible cobordism. 
 
We first build a particular extension of $\Lambda$ to the vertical 
3-skeleton of 
$M\times I$ this way. Given any  vertical 3-cube remark that $\Lambda$ 
evaluates non-trivially on an even number of 
its 4 vertical 2-faces (possibly none). Trace a vertical 8 on any 
vertical face with non-trivial $\Lambda$, and connect the 8's in pairs 
by means of tubes whose section is a vertical 8. 
If all of the 4 faces are traced, the pairs must be of adjacent 
faces. 
Now remark that this 
first extension is further extendible to the vertical 4-skeleton. 
Indeed consider any vertical 4-cube. The collection of its 6 vertical 
3-faces (that can be visualized as an $S^2\times (0,1)$) 
contains a disjoint union of immersions each representing  an element of 
$2P_3$. But each of these immersions is the trivial 
element. Indeed 
by construction the top of the 8 describes a curve which bounds a disk not 
intersecting the curve of double points, 
 hence 
having trivial linking number with it, and the double of this linking number is the Arf invariant of the immersion. 
Hence 
$\Lambda$ is extendible to the vertical 4-skeleton. The following 
obstructions are all trivial by theorem~\ref{forte}, so $\Lambda$ is 
extendible to $M\times (0,1)$, as it was claimed. 
 
Now given a $s$-admissible cobordism $g$ between 
$(s+1)$-admissible immersions, $3<s\leq7$, define the same way the 
$s$-cochain $\Lambda$ closed in $M\times (0,1)$ and by directly 
applying theorem~\ref{forte} extend it to an immersion $h$ in 
$M\times(0,1)$. Up to adding $h$ an appropriate number of times one 
obtains a cobordism $g'$ that is still $s$-admissible but such that 
$g'(e^{s-1})=0$ for all $(s-1)$-cube $e^{s-1}$ of $M$, hence that can be deformed to 
a $(s+1)$-admissible cobordism.~$\square$ 
 
\end{dimostraz}

\section{The group structure on orientable 4-manifolds}

The graded group $gr^*(N(M))$ is isomorphic to $N(M)$ as a set, but looses 
its group structure. We give a result concerning the group structure 
when $M$ is an orientable 4-manifold. 
 
We first remark  that for any $n$-manifold the 
{\em total James-Hopf invariant} $JH$, 
that is, the product of the  James-Hopf invariants composed with 
Poincar\'e duality, becomes a homomorphism of groups with 
$\oplus_{j=1}^n 
H^j(M,P_j)$ endowed with  the group structure coming from that of 
algebra 
\[ 
(\alpha*\beta)_j=\alpha_j+\beta_j+\sum_{s+t=j}\alpha_s\beta_t. 
\] 
In dimension 3 the invariant $JH$ provides completely the group 
structure, up to immersions that are contained in a ball. These 
last form a subgroup that can be detected by a version of the 
Arf invariant of $P_3$ (though Benedetti and Silhol provided 
a deeper invariant). 
Up to immersions in a ball, any class can be realized as the decoration of an embedded 
representative and any immersion can be split in a unique way as the 
sum of immersions obtained by decorating a submanifold, that is, $JH$ 
is surjective and injective. Neither of those properties hold for 
orientable 4-manifolds. Indeed from the main theorem self-intersection 
of 2-classes is the (only) obstruction for decorating an embedded 
representative. Moreover decorating a simple curve with an element of 
$2P_3$ provides an immersion in $F^3$ with trivial $JH$ but non-trivial 
$\chi^3$. 
 
We define the map $k$ from $H_1(M,\Z/4\Z)=H_1(M,\Z)\otimes \Z/4\Z$ in $N(M)$ 
that associates to $\gamma\otimes m$ an embedded curve representing $\gamma$ decorated by 
the canonical 
immersion with invariant $2m$. This map results to be 
well-defined, that is, $\gamma\ltimes 2m$ does not depend either on the 
representative of $\gamma$ nor on the trivialization of the normal 
bundle. 
Moreover the image of $JH$ is the subgroup of $\oplus_{j=1}^4 
H^j(M,P_j)$ whose support is $H^1(M,\zd)\oplus Q^2(M,\zd)\oplus 
H^3(M,\zd)\oplus H^4(M,\zd)$ (see~\cite{GinPhD} for more details). The 
following holds. 
 
\begin{prop} There is a short 
exact sequence of groups 
\[ 
0\rightarrow H_1(M,\Z/4\Z)\stackrel{k}{\rightarrow} 
N(M)\stackrel{JH}{\rightarrow}H^1(M,\zd)* Q^2(M,\zd)* 
H^3(M,\zd)* H^4(M,\zd)\rightarrow 0. 
\] 
\end{prop} 
 
\begin{dimostraz} That $k$ is injective follows from the fact that 
images of different cycles have different $\chi^3$ and this last is 
injective. 
We show exactness in the middle term. That $JH(k(\alpha\otimes m))$ is trivial 
for any $\alpha\otimes n\in H_1(M,\Z)\otimes\Z/4\Z$ follows 
since any representative $f$ of $k(\alpha\otimes m)$ 
 has no triple points nor quadruple points, and the locus of its double points 
is a surface representing the trivial element of $H_2(M,\zd)$. The whole of $f$ retracts on the decorated curve in fact, 
hence also $JH_3(f)$ is trivial. 
Now suppose that $JH(f)=(0,0,0,0)$.  We must prove 
that $f$ is in the image of the map $k$. 
Since $\chi^1(f)=JH_3(f)=0$, from lemma~\ref{nucleo} $f$ 
belongs to $F_2$. Assume then that $f$ is 2-admissible. 
Now $f$ has trivial second cohomological 
invariant, since $\chi^2(f)=PD(JH_2(f))$, 
hence in particular $f$ belongs to $F^3$, see 
proposition~\ref{nucleo}. Assume that $f$ is 3-admissible and 
consider 
$\chi^3(f)\in H^3(M,P_3)$. This class has trivial 
reduction modulo 2, since 
$$ 
\chi^3(f)(mod\; 2)=PD(JH_1(f))=0, 
$$ 
so there is an 
element  $\kappa_f\in H_1(M,\Z/4\Z)$ such that 
$\chi^3(f)=2PD(\kappa_f)$. It is easy to see that 
$f=k(\kappa_f)$.~$\square$ 
 
\end{dimostraz} 
 
\begin{ex}{\rm 
The group of the complex projective plane is 
$N({\bf CP}^2)=\zd$ generated by a non-trivial immersion in a small 
ball. 
} 
\end{ex} 
 
\begin{oss}{\rm 
The cobordism group of a 4-manifold is generated by embedded decorated 
submanifolds. It is proven in \cite{BeSSPS} that the cobordism 
class of a codimension-one embedding only depends on homology 
modulo 2, and in \cite{GinPhD} that the cobordism class of a 
decorated curve only depends on the (oriented) homology of the curve and the 
cobordism class of the decorating immersion. So if $M$ is an oriented 4-manifold such that 
$Q^2(M,\zd)=0$ one can choose a the set of 
generators of $N(M)$ by choosing a set $\{C_j\}$ of oriented 
curves generating $H_1(M,\Z)$ and a set $\{S_k\}$ of 
codimension-one embeddings generating $H^1(M,\zd)$ and considering 
the set $\{\cdot\ltimes P, C_j\ltimes B, S_k\ltimes\cdot\}$, $B$ denoting the left 
Boy immersion and $P$ denoting a  generator of $P_4$.

 
 
} 
\end{oss} 
 
\begin{oss} {\rm The image of $PD(JH_2)$ being $Q^2(M,\zd)$ implies in 
particular that for any pair of 1-cocycles modulo 2, $\alpha$ and 
$\beta$, the relation $\alpha^2\beta^2=0$ holds. This fact has an 
elementary 
geometric proof. Represent the dual of $\alpha$ and $\beta$ by means of 
embedded hypersurfaces $A$ and $B$. Remark that if $A$ is orientable 
then $\alpha^2$ is 0, and if $A$ is non-orientable its 
self-intersection is the orientation cycle of $A$, hence can be 
represented by an orientable surface $F$ in $A$. Call $C$ the curve 
intersection between $F$ and $B$. 
Then a 
representative of the dual of $\alpha^2\beta^2$ is the intersection 
between $C$ and $B$. But the normal bundle to $B$ restricted to $C$ 
is orientable since it is the normal bundle to $C$ in $F$, that is 
trivial since $F$ is orientable, and hence $C\cdot B=0$. 
} 
\end{oss}

 \begin{small} 
\bibliographystyle{plain}

\begin{thebibliography}{24} 
 
\bibitem{BeSSPS} 
R.~Benedetti and R.~Silhol, 
\newblock {\em ${S}pin$ and ${P}in^-$ structures, immersed and embedded surfaces and 
  a result of {S}egre on real cubic surfaces}, 
\newblock {Topology}, {\bf 34}(1995), 651--678. 
 
\bibitem{CarGBS} 
J.~S. Carter, 
\newblock {\em On generalizing {B}oy's surface: constructing a generator of the 
  third stable stem}, 
\newblock {Trans.A.M.S.}, {\bf 298}(1986), 103--122. 
 
\bibitem{FunCIM} 
L.~Funar, 
\newblock {\em Cubulations, immersions, mappability and a problem of 
{H}abegger}, 
\newblock {Ann. Sci. {\'E}cole Norm. Sup}, {\bf 32}(1999), 681--700. 
 
 
 
\bibitem{GinPhD} 
R.~Gini, 
\newblock {\em Cobordism of codimension-one immersions}, 
\newblock PhD.~Thesis, Universit\`a di Pisa, Dipartimento di 
Matematica, 2001. 
 
\bibitem{rosa} 
R.~Gini, 
\newblock {\em Cobordism of immersions of surfaces in non-orientable 3-manifolds}, 
\newblock { Manuscripta Math.}, {\bf 104}(2001), 49-69. 
 
\bibitem{HilGCT} 
P.~Hilton, 
\newblock {General Cohomology Theory and {K}-Theory}, 
\newblock v.1, {\em Lecture Note Series, Cambridge University Press}, 
  1971. 
 
\bibitem{JamOST} 
I.~M. James, 
\newblock {\em On the suspension triad}, 
\newblock {Ann. of Math.}, {\bf 63}(1956), 191--247. 
 
 
\bibitem{HaHIST} 
J.~Hass and J.~Hughes, 
\newblock {\em Immersions of surfaces in 3-manifolds}, 
\newblock {Topology}, {\bf 24}(1985), 97-112. 
 
 
\bibitem{KirPLT} 
R.~Kirby, 
\newblock Problems in low-dimensional topology, 
\newblock {\em {"}Geometric Topology{"}, Georgia 
  International Topology Conference}, 
 (W.H.Kazez, Ed.), 
 vol.~2, {\em Studies in Advanced 
  Math}, p. 35--472. AMS-IP, 1995. 
 
\bibitem{KoSSIH} 
U.~Koschorke and B.~Sanderson, 
\newblock {\em Self--intersections and higher {H}opf invariants}, 
\newblock {Topology}, {\bf 17}(1978), 283--290. 
 
\bibitem{LiuTHA} 
A.~Liulevicius, 
\newblock {\em A theorem in homological algebra and stable homotopy of projective 
  spaces}, 
\newblock {Trans.A.M.S.}, {\bf 109}(1963), 540--552. 
 
\bibitem{MasSHT} 
W.~S. Massey, 
\newblock {Singular Homology Theory}, 
\newblock {\em Graduate Texts in Math., Springer},  1980. 
 
 
 
\bibitem{NowQPR} 
T. ~Nowik, 
\newblock {\em Quadruple Points of Regular Homotopies of Surfaces in 
3-Manifolds}, 
\newblock   {Topology}, {\bf 39}(2000), 1069-1088. 
 
 
 
\bibitem{PicLHT} 
R.~A. Piccinini, 
\newblock {Lectures on Homotopy Theory}, 
\newblock {\em North-Holland},  1992. 
 
\bibitem{PinRHC} 
U.~Pinkall, 
\newblock {\em Regular homotopy classes of immersed surface}, 
\newblock {Topology}, {\bf 24}(1985), 421--434. 
 
\bibitem{ThoQPG} 
R.~Thom, 
\newblock {\em Quelques propri\'et\'es globales des vari\'et\'es differentiables}, 
\newblock {Comment. Math. Helvetici}, {\bf 28}(1954), 17--86. 
 
\bibitem{VogCI} 
P.~Vogel, 
\newblock {\em Cobordisme d'immersions}, 
\newblock {Ann.Sci.\'Ecole Norm.Sup.},  {\bf 7}(1974), 317-358. 
 
 
 
\bibitem{WelCGI} 
R.~Wells, 
\newblock {\em Cobordism groups of immersions}, 
\newblock {Topology}, {\bf 5}(1966), 281-294. 
 
\bibitem{WhiEHT} 
G.~W. Whitehead, 
\newblock {Elements of Homotopy Theory}, 
\newblock {\em Graduate Texts in Math., Springer}, 1978. 
 
 
 
\end{thebibliography}

\end{small}

 \end{document}